\documentclass[11pt]{article}
\usepackage[utf8]{inputenc}
\usepackage{graphicx}
\usepackage{caption}
\usepackage{amsmath}
\usepackage{amssymb}
\usepackage{float}
\usepackage{amsthm}
\usepackage{xcolor}
\usepackage{comment}

\usepackage{mathtools} %%%%%%%%% num?rote uniquement les ?quations r?f?renc?es
\mathtoolsset{showonlyrefs} %%%%%%%%%%%%%%%

\DeclareMathOperator{\Tr}{{\text{Tr}}}

\DeclareMathOperator{\R}{\mathbb{R}}

\def\Ec{{\cal E}}
\def\Lc{{\cal L}}
\def\Nc{{\cal N}}

\def\Fc{{\cal F}}
\def\Uc{{\cal U}}
\def\Vc{{\cal V}}
\def\Wc{{\cal W}}
\def\Xc{{\cal X}}

\def\Zc{{\cal Z}}

\def \E{\mathbb{E}}
\def \F{\mathbb{F}}
\def \M{\mathbb{M}}
\def \P{\mathbb{P}}

\def\1{{\bf 1}}
\def \I{{\bf I}}
\def \N{\mathbb{N}}

\newcommand*\diff{\mathop{}\!\mathrm{d}}

\def \eps{\varepsilon}

\def\argmin{\mathop{\rm argmin}}
\def\argmin_#1{\underset{#1}{\mathrm{argmin\, }}}

\def \ep{\hbox{ }\hfill$\Box$}

\newcommand{\un}{1\hspace{-1mm}{\rm I}}   % 1

\newtheorem{Theorem}{Theorem}[section]

\newtheorem{Remark}{Remark}[section]

\numberwithin{equation}{section} % Num?rotation des ?quations en rapport avec la section.

\usepackage[backend=bibtex, style=alphabetic]{biblatex}
\usepackage{hyperref}
\def \trans{^{\scriptscriptstyle{\intercal}}}

\usepackage{lineno} 
\begin{document}
\title{Deep backward schemes for \\ high-dimensional nonlinear PDEs
\thanks{This work is supported by  FiME, Laboratoire de Finance des March\'es de l'Energie, and the ''Finance and Sustainable Development'' EDF - CACIB Chair.}}

\author{C\^ome \sc{Hur\'e}\footnote{LPSM, Paris-Diderot University \sf \href{mailto:hure at lpsm.paris}{hure at lpsm.paris}}  \and
Huy\^en \sc{Pham}\footnote{LPSM, Paris-Diderot University, CREST-ENSAE \& FiME \sf \href{mailto:pham at lpsm.paris}{pham at lpsm.paris}} 
\and
Xavier \sc{Warin}
\footnote{EDF R\&D \& FiME \sf \href{mailto:xavier.warin at edf.fr}{xavier.warin at edf.fr}} 
}

\date{First version: february 4, 2019 \\
This final accepted version: November 8, 2019}

\maketitle              % typeset the header of the contribution

\begin{abstract}
We propose new machine learning schemes  for solving high dimensional nonlinear partial differential equations (PDEs). Relying on the classical backward stochastic di\-fferential equation (BSDE) representation of PDEs, our algorithms estimate simultaneously the solution and its gradient by deep neural networks. 
%In contrast with the recent app\-roach in \cite{weinan2017deep}, 
These approximations are performed at each time step  from the minimization of loss functions defined recursively by backward induction. The metho\-dology is extended to variational inequalities ari\-sing in optimal stopping problems. 
We analyze the convergence of the deep learning schemes and provide error estimates in terms of the universal approximation of neural networks.  
Numerical results  show that our algorithms give very good results till dimension 50 (and certainly above), for both PDEs and  variational inequalities problems.  
%that the recent method in \cite{weinan2017deep} cannot solve.
For the PDEs resolution, our results are very similar to those obtained by the recent method in \cite{weinan2017deep} when the latter converges to the right solution or does not diverge. Numerical tests indicate that the proposed methods are not 
stuck in poor local minima
%trapped by local minima far away from the true solution 
as it can be the case with the algorithm designed in \cite{weinan2017deep}, and no divergence is experienced. The only limitation seems to be due to the inability of the considered deep neural networks to represent a solution with a too complex structure in high dimension.
\end{abstract}

\vspace{5mm}

\noindent {\bf Key words:} Deep neural networks,  nonlinear PDEs in high dimension, optimal sto\-pping problem, backward stochastic differential equations.  

\vspace{5mm}

\noindent {\bf MSC Classification:}  60H35, 65C20, 65M12.

\newpage

%\renewcommand{\abstractname}{Acknowledgements}
%\begin{abstract}
%The author(s) acknowledge support from the Chair Finance \& Sustainable Development and/or the FiME Lab  (Institut Europlace de Finance)
%\end{abstract}
%
%
%

\section{Introduction}
This paper is devoted to the resolution in high dimension of nonlinear parabolic partial differential equations (PDEs) of the form
\begin{equation} \label{eq:PDEInit}
\left\{
\begin{aligned}
 \partial_t u + \Lc u + f(.,.,u,\sigma\trans D_x u) & = 0 , \;\;\;\;\;\; \mbox{ on } [ 0,T)\times\R^d, \\
 u(T,.) &=g, \;\;\;\;\;  \mbox{ on } \R^d,
 \end{aligned}
 \right.
 \end{equation}
with a  non-linearity  in the solution and its gradient via the function $f(t,x,y,z)$ defined on $[0,T]\times\R^d\times\R\times\R^d$, a  terminal condition $g$,   
and a second-order generator $\Lc$ defined by
 \begin{align}
 \label{eq:PDE}
 \mathcal{L}u   & := \frac{1}{2} \Tr \big(\sigma \sigma\trans D_x^2 u \big) + \mu.D_x u.
 \end{align}
Here $\mu$ is a function defined on $[0,T] \times \R^d$ with values in $\R^d$, $\sigma$ is a function defined on $[0,T] \times \R^d$ with values in  $\M^d$ the set of  $d \times d$ matrices, and  $\mathcal{L}$ is the generator associated to the forward diffusion process:
  \begin{flalign}
  \Xc_t = x_0 + \int_0^t \mu(s,\Xc_s) \diff s+  \int_0^t  \sigma(s,\Xc_s)  \diff W_s, \;\;\; 0 \leq t \leq T, 
 \label{eq:SDE}
 \end{flalign}
 with  $W$ a $d$-dimensional Brownian motion on some probability space $(\Omega,\Fc,\P)$ equipped with a filtration 
 $\F$ $=$ $(\Fc_t)_{0\leq t\leq T}$ satisfying the usual conditions. 
 
 Due to the so called ``curse of dimensionality", the resolution of nonlinear PDEs in high dimension has always been a challenge for scientists.
Until recently, only the BSDE (Backward Stochastic Differential Equation) approach first developed in  \cite{pardoux1990adapted} was avai\-lable to tackle this problem: using the time discretization scheme proposed in \cite{bouchard2004discrete}, some effective algorithms based on regressions manage to solve non linear PDEs in dimension above 4 (see \cite{gobet2005regression,lemor2006rate}).
However this approach is still not implementable in dimension above 6 or 7 : the number of basis functions used for the regression still explodes with the dimension.

Quite recently some new methods have been developed for this problem, and  several methodologies have emerged:
\begin{itemize}
    \item Some are based on the Feyman-Kac representation of the PDE. Branching techniques \cite{henry2016branching} have been studied and shown to be convergent but only for small maturities and some small nonlinearities. Some  effective techniques based on nes\-ting Monte Carlo have been studied in \cite{warin2018nestingsMC,warin2018monte}: the convergence is proved for semi-linear equations.
    Still based on this Feyman-Kac representation some machine learning techniques  permitting to solve a fixed point problem have been used recently in \cite{chan2018machine}: numerical results show that it is efficient and some partial demonstrations justify why it is effective.
    \item Multilevel Picard methods have been developed in  \cite{hutetal18} and \cite{hutzenthaler2018overcoming} with algorithms based on Picard iterations, multi-level techniques and automatic differentiation. These methods permit to handle some high dimensional PDEs with non linearity in $u$ and its gradient $D_x u$, with  convergence results as well as numerous numerical examples showing their efficiency in high dimension. 
    \item Another class of methods is based on the BSDE approach and the curse of dimensionality issue is partially avoided by using some machine learning techniques. 
    The pioneering papers \cite{han2017overcoming,weinan2017deep} propose a neural-networks based technique called \emph{Deep BSDE}, which was  the first serious attempt for using machine learning methods to solve high dimensional PDEs.  Based on an Euler discretization of the forward underlying SDE  $\Xc_t$, the idea is to view the BSDE as a forward SDE, and the algorithm tries to learn the values $u$ and  $z$ $=$ $\sigma\trans Du$ at each time step of the Euler scheme by minimizing a global loss function between the forward simulation of $u$ till maturity $T$ and  the target $g(\Xc_T)$. This deep learning approximation has been extended to the case of fully nonlinear PDE and second order BSDE in \cite{becketal19}.
    \item At last, using some machine learning representation of the solution,  \cite{sirignano2018dgm} proposes with the so-called Deep Galerkin Method  to use the automatic numerical differentiation of the solution to solve the PDE on a finite domain. 
   The authors prove the convergence of their method but without information on the rate of convergence. 
\end{itemize}

Like the second methodology, our approach  relies on BSDE representation of the PDE and deep learning appro\-ximations: we first discretize the BSDE associated to the PDE by an Euler scheme, but in contrast with \cite{weinan2017deep}, we  adopt a classical backward resolution technique. 
On each time step,  we propose to use  some machine learning techniques  to estimate simultaneously the solution and its gradient by minimizing a loss function defined recursively by backward induction,  and solving this local problem by a stochastic gradient algorithm. 
Two different schemes are designed to deal with  the local problems:
\begin{itemize}
    \item[(1)] The first one tries the estimate the solution and its gradient by a neural network.
    \item[(2)] The second one tries only to approximate the solution by a neural network while its gradient is estimated directly with some numerical differentiation techniques.
\end{itemize}
The proposed methodology is then extended to solve some variational inequalities, i.e., free boundary problems related to optimal stopping problems. We mention that the related recent paper \cite{becker} also proposes deep learning method for solving optimal stopping problems, but differently from our method, it relies on the approximation of (randomised) stopping decisions with a sequence of multilayer feedforward neural networks.

Convergence analysis of the two schemes for PDEs and variational inequalities is provided and shows that the  approximation error goes to zero as we increase the number of time steps and the number of neurons/layers whenever the gradient descent method used to solve the local problems is not trapped in a local minimum.  
Notice that similar convergence result for the deep BSDE method has been also obtained in \cite{hanlong18} with a posteriori error estimation of the solution in terms of the  universal approximation capability of global neural networks.

In the last part of the paper, we test our algorithms on different examples. When the solution is easy to represent by a neural network,  we can solve the problem in quite  high dimension (at least $50$ in our numerical tests). 
We show that the proposed methodology improves the algorithm proposed in \cite{han2017overcoming} that sometimes does not converge or is trapped in a local minimum far away from the true solution.
We then show that when the solution has a very complex structure, we can still solve the problem but only in moderate dimension: the neural network used is  not anymore able to represent the solution accurately in very  high dimension. 
Finally, we illustrate numerically that the method is effective to solve some system of varia\-tional inequalities: we consider the problem of American options and show that it can be solved very accurately in high dimension (we tested until $40$). 

The outline of the paper is organized as follows. In Section \ref{secNN}, we give a brief and useful reminder for neural networks. We describe in Section \ref{secalgo} our two numerical schemes and compare with the algorithm in 
\cite{han2017overcoming}. Section \ref{secconv} is devoted to the convergence analysis of our machine learning algorithms, and we present in Section \ref{secnum} several numerical tests.

\section{Neural networks as function approximators}
\label{secNN}

Multilayer (also called deep) neural  networks are designed to approximate unknown or large class of functions. In contrast to additive approximation theory with  weighted sum over basis functions, e.g. polynomials,  
neural networks  rely on the composition of simple functions, and appear to provide an efficient way to handle high-dimensional approximation problems, in particular thanks to the increase in computer power for finding the ``optimal" parameters by (stochastic) gradient descent methods.      

We shall consider feedforward (or artificial) neural networks, which represent  the basic type of deep neural networks. Let us recall some notation and basic definitions that will be useful in our context.  We fix the input dimension $d_0$ $=$ $d$ (here the dimension of the state variable $x$), the output dimension $d_1$ (here $d_1$ $=$ $1$  for approximating the real-valued solution to the PDE, or $d_1$ $=$ $d$ for approximating the vector-valued gradient function), the global number $L+1$ $\in$ $\N\setminus\{1,2\}$ of layers with $m_\ell$, $\ell$ $=$ $0,\ldots,L$, the number of neurons (units or nodes) on each layer: the first layer is the input layer with $m_0$ $=$ $d$, the last layer is the output layer with $m_L$ $=$ $d_1$, and the $L-1$ layers between are called hidden layers, where we choose for simplicity the same  dimension $m_\ell$ $=$ $m$, $\ell$ $=$ $1,\ldots,L-1$.

A  feedforward neural network is a  function from  $\R^{d}$ to $\R^{d_1}$ defined as the composition
\begin{align} \label{defNN}
x \in \R^d  & \longmapsto  \; A_L \circ  \varrho \circ A_{L - 1} \circ \ldots \circ \varrho \circ A_1(x) \; \in \; \R^{d_1}. 
\end{align}
Here $A_\ell$, $\ell$ $=$ $1,\ldots,L$ are affine transformations: $A_1$ maps from $\R^d$ to $\R^m$, $A_2,\ldots,A_{L-1}$ map from $\R^m$ to $\R^m$, and $A_L$ maps from $\R^m$ to 
$\R^{d_1}$, represented by 
\begin{align}
A_\ell (x) &= \; \Wc_\ell x + \beta_\ell,
\end{align}
for a matrix $\Wc_\ell$ called weight, and a vector $\beta_\ell$ called  bias term,  $\varrho$ $:$ $\R$ $\rightarrow$ $\R$ is a nonlinear function, called activation function, and applied 
component-wise on the outputs of $A_\ell$, i.e., $\varrho(x_1,\ldots,x_m)$ $=$ $(\varrho(x_1),\ldots,\varrho(x_m))$. Standard examples of activation functions are 
the sigmoid, the ReLu, the Elu, $\tanh$.

All these matrices $\Wc_\ell$ and vectors $\beta_\ell$, $\ell$ $=$ 
$1,\ldots,L$,  are the parameters of the neural network, and can be identified with  an element $\theta$ $\in$ $\R^{N_m}$, where $N_m$ $=$ 
$\sum_{\ell=0}^{L-1} m_\ell (1+m_{\ell+1})$ $=$ $d(1+m)+m(1+m)(L-2)+m(1+d_1)$ is the number of parameters, where we fix $d_0$, $d_1$, $L$, but allow  
growing number $m$ of hidden neurons.  
We denote by $\Theta_m$ the set of possible parameters: in the sequel, we shall consider either the case when there are no constraints on parameters, i.e.,  
$\Theta_m$ $=$ $\R^{N_m}$, or when the total variation norm of the neural networks is smaller than $\gamma_m$, i.e., 
\begin{align}
\Theta_m & = \Theta_m^\gamma  \; := \;  
 \big\{ \theta = (\Wc_\ell,\beta_\ell)_\ell:  |\Wc_l| \leq \gamma_m, \;\; \ell = 1,\ldots,L \big\}, \;\; \mbox{ with} \; \gamma_m \nearrow \infty, \mbox{ as } m \rightarrow \infty. 
\end{align}
We denote by $\Phi_{_m}(.;\theta)$ the neural network function defined in \eqref{defNN}, and  by $\Nc\Nc_{d,d_1,L,m}^\varrho(\Theta_m)$ the set of all such neural networks $\Phi_{_m}(.;\theta)$ 
for $\theta$ $\in$ $\Theta_m$, and set 
\begin{eqnarray*}
\Nc\Nc_{d,d_1,L}^\varrho &= &  
\bigcup_{m \in \N}  \Nc\Nc_{d,d_1,L,m}^\varrho(\Theta_m) \; = \; \bigcup_{m \in \N}  \Nc\Nc_{d,d_1,L,m}^\varrho(\R^{N_m}),     
\end{eqnarray*}
as  the class of all neural networks within a fixed structure given by $d$,  $d_1$, $L$ and $\varrho$. 
 
The fundamental result of Hornick et al. \cite{horetal89} justifies the use of neural networks as function approximators:
 
 \vspace{1mm}
 
\noindent  {\bf Universal approximation theorem (I)}:  $\Nc\Nc_{d,d_1,L}^\varrho$ is dense in $L^2(\nu)$ for any finite measure $\nu$ on $\R^d$, whenever $\varrho$ is continuous and non-constant. 

\vspace{1mm}

Moreover, we  have a universal approximation result for the derivatives in the case of a single hidden layer, i.e. $L$ $=$ $2$, and when the activation function is a smooth function, see \cite{hor90universal}. 

\vspace{1mm}

\noindent  {\bf Universal approximation theorem (II)}:  Assume that $\varrho$ is a (non constant) $C^k$ function. Then, $\Nc\Nc_{d,d_1,2}^\varrho$ approximates any function and its derivatives up to order $k$,  arbitrary well on any compact set of $\R^d$.

\section{Deep learning-based schemes for semi-linear PDEs}
\label{secalgo}

The starting point for our probabilistic numerical schemes to the PDE \eqref{eq:PDEInit} is the well-known (see \cite{pardoux1990adapted}) nonlinear Feynman-Kac formula via the pair $(Y,Z)$ of 
$\F$-adapted processes valued in  $\R\times\R^d$,  solution to the BSDE  
\begin{align} \label{eqBSDE}
Y_t   &= g(\Xc_T) + \int_t^T f(s,\Xc_s,Y_s,Z_s) \diff s - \int_t^T Z_s\trans  \diff W_s, \;\;\; 0 \leq t \leq T,
\end{align}
related to the solution $u$ of \eqref{eq:PDEInit} via
\begin{align*}
Y_t  & = u(t,\Xc_t), \;\;\; 0 \leq t \leq T,
\end{align*}
and when $u$ is smooth:
\begin{align*}
Z_t & = \;  \sigma\trans(t,\Xc_t)D_x u(t,\Xc_t), \;\;\; 0 \leq t \leq T. 
\end{align*}

\subsection{The deep BSDE scheme of \cite{han2017overcoming}}

The DBSDE algorithm proposed in \cite{han2017overcoming,weinan2017deep} starts from the BSDE representation \eqref{eqBSDE} of the solution to  \eqref{eq:PDEInit}, but rewritten in forward form as: 
\begin{align}
\label{eq:bsde}
u(t,\Xc_t) = & \;\; u(0, x_0)- \int_0^t f(s,\Xc_s,u(s,\Xc_s),\sigma\trans(s,\Xc_s)D_xu(s,\Xc_s)) \diff s \\
& \;\;\; + \int_0^t D_xu(s,\Xc_s)\trans\sigma(s,\Xc_s) \diff W_s, 
\;\;\;\;\; 0 \leq t \leq T. 
\end{align}

The forward process $\Xc$ in equation \eqref{eq:SDE}, when it is not simulatable, is numerically approximated by an Euler scheme $X$ $=$ $X^\pi$ 
on a time grid: $\pi$ $=$ $\{t_0=0<t_1< \ldots < t_N = T\}$, with modulus $|\pi|$ $=$ $\max_{i=0,\ldots,N-1}\Delta t_i$, $\Delta t_i$ $:=$ $t_{i+1}-t_i$, and defined as 
%(we keep the same notation $X$ for the Euler scheme)  
\begin{equation}
X_{t_{i+1}} \; = \;   X_{t_{i}} + \mu(t_i, X_{t_i}) \Delta t_{i} + \sigma(t_i, X_{t_i}) \Delta W_{t_i}, \;\;\; i=0,\ldots,N-1, \; X_0 = x_0, 
\label{eq:eulerSDE}
\end{equation}
where we set $\Delta W_{t_i}$ $:=$ $W_{t_{i+1}} - W_{t_{i}}$.  To alleviate notations, we omit the dependence of $X$ $=$ $X^\pi$ on the time grid $\pi$ as there is no ambiguity 
(recall that we use the notation $\Xc$ for the  forward diffusion process). 
The approximation of equation \eqref{eq:PDEInit} is then given formally from the Euler scheme associated to the forward representation \eqref{eq:bsde} by 
\begin{align} \label{uF}
u(t_{i+1}, X_{t_{i+1}}) & \approx  F(t_i,X_{t_i},u(t_i,X_{t_i}),\sigma\trans(t_i,X_{t_i})D_x u(t_{i}, X_{t_{i}}),\Delta t_{i},\Delta W_{t_i})
\end{align}
with 
\begin{align} 
F(t,x,y,z,h,\Delta) & :=  y - f(t,x,y,z)h + z\trans\Delta.   
\label{defF}
\end{align}
%In the machine learning language, the independent realizations of $(X_{t_i})_{i =1...N}$ represent the training %data.

In \cite{han2017overcoming,weinan2017deep}, the numerical approximation of $u(t_i,X_{t_i})$ is designed as follows: 
starting from an estimation $\Uc_0$ of $u(0,X_0)$, and then using at each time step $t_i$, $i$ $=$ $0,\ldots,N-1$, a  multilayer neural network $x$ $\in$ $\R^d$ $\mapsto$ $\Zc_i(x;\theta_i)$ with parameter $\theta_i$ for the approximation of $x$ $\mapsto$ 
$\sigma\trans(t_i,x)D_x u(t_i,x)$: 
\begin{align} \label{NNti} 
\Zc_i(x;\theta_i) & \approx \sigma\trans(t_i,x)D_x u(t_i,x), 
\end{align}
one computes estimations $\Uc_i$ of $u(t_i;X_{t_i})$ by forward induction via:
\begin{align*}
\Uc_{i+1} &=  F(t_i,X_{t_i},\Uc_i,\Zc_i(X_{t_i};\theta_i),\Delta t_i,\Delta W_{t_i}),   
\end{align*}
for $i$ $=$ $0,\ldots,N-1$. This algorithm forms a global deep neural network composed of the neural networks \eqref{NNti} of each period, by taking as input data (in machine learning language)  the paths of $(X_{t_i})_{i=0,\ldots,N}$ and $(W_{t_i})_{i=0,\ldots,N}$, and giving as output $\Uc_N$ $=$ $\Uc_N(\theta)$, which is a function  of the input and of the total set of parameters $\theta$ $=$ $(\Uc_0,\theta_0,\ldots,\theta_{N-1})$.  The output aims to match the terminal condition $g(X_{t_N})$ of the BSDE, and one then optimizes over the parameter $\theta$ the expected square loss function:
\begin{align*}
\theta & \mapsto \; \E \big| g(X_{t_N})  - \Uc_N(\theta) \big|^2. 
\end{align*}
This is obtained by stochastic gradient descent-type (SGD) algorithms relying on training input data.

\subsection{New schemes: DBDP1 and DBDP2} \label{secnewML}

The proposed scheme is defined from a backward dynamic programming type relation, and  has two versions:

\begin{itemize}
    \item[(1)] First version: 
    \begin{itemize}
    	\item[-] Initialize from an estimation $\widehat\Uc_N^{(1)}$ of $u(t_N,.)$ with $\widehat\Uc_N^{(1)}$ $=$ $g$
    	\item[-] For $i$ $=$ $N-1,\ldots,0$, given $\widehat\Uc_{i+1}^{(1)}$, use a pair of deep neural networks $(\Uc_i(.;\theta),\Zc_i(.;\theta))$ $\in$ $\Nc\Nc_{d,1,L,m}^\varrho(\R^{N_m})\times\Nc\Nc_{d,d,L,m}^\varrho(\R^{N_m})$  for the approximation of $(u(t_i,.),\sigma\trans(t_i,.)D_x u(t_i,.))$, and compute (by SGD) the minimizer of the expected quadratic loss function
    \begin{equation} \label{eq:scheme1}
    \left\{ 
    \begin{aligned} 
    \hat L_i^{(1)}(\theta) & := \E \Big| \widehat\Uc_{i+1}^{(1)}(X_{t_{i+1}}) - F(t_i,X_{t_i},\Uc_i(X_{t_i};\theta),\Zc_i(X_{t_i};\theta),\Delta t_i,\Delta W_{t_i}) \Big|^2  \\
    \theta_i^* & \in {\rm arg}\min_{\theta\in\R^{N_m}} \hat L_i^1(\theta). 
    \end{aligned}
    \right.
    \end{equation}
      Then, update: $\widehat\Uc_i^{(1)}$ $=$ $\Uc_i(.;\theta_i^*)$, and set $\widehat\Zc_i^{(1)}$ $=$ $\Zc_i(.;\theta_i^*)$. 
    \end{itemize}
     \item[(2)] Second version: 
     \begin{itemize}
     \item Initialize with $\widehat\Uc_N^{(2)}$ $=$ $g$ 
     \item For $i$ $=$ $N-1,\ldots,0$, given $\widehat\Uc_{i+1}^{(2)}$, use a deep neural network $\Uc_i(.;\theta)$ $\in$ $\Nc\Nc_{d,1,L,m}^\varrho(\Theta_m)$, 
     and compute (by SGD) the minimizer of the expected quadratic loss function
     \begin{equation} \label{eq:scheme2}
     \left\{ 
     \begin{aligned} 
     \hat L_i^{(2)}(\theta) & := \E \Big| \widehat\Uc_{i+1}^{(2)}(X_{t_{i+1}}) -  \\
     & \quad \quad F(t_i,X_{t_i},\Uc_i(X_{t_i};\theta),\sigma\trans(t_i,X_{t_i}) \hat D_x \Uc_i(X_{t_i};\theta),\Delta t_i,\Delta W_{t_i}) \Big|^2 \\
     \theta_i^* & \in {\rm arg}\min_{\theta\in\Theta_m} \hat L_i^2(\theta), 
     \end{aligned}
     \right.
     \end{equation}
     where $\hat D_x \Uc_i(.;\theta)$ is the numerical differentiation of $\Uc_i(.;\theta)$.  
     Then, update: $\widehat\Uc_i^{(2)}$ $=$ $\Uc_i(.;\theta_i^*)$, and set $\widehat\Zc_i^{(2)}$ $=$ $\sigma\trans(t_i,.) \hat D_x \Uc_i(.;\theta_i^*)$. 
     \end{itemize}
\end{itemize}

\begin{Remark} \label{remNN}
{\rm For the first version of the scheme, one can use independent neural networks, respectively for the approximation of $u(t_i,.)$ and for the approximation of $\sigma\trans(t_i,.)D_xu(t_i,.)$. 
In other words,  the parameters are divided into a pair $\theta$ $=$ $(\xi,\eta)$ and we consider neural networks $\Uc_i(.;\xi)$ and $\Zc_i(.;\eta)$. 
}
\ep
\end{Remark}
In the sequel, we refer to the first and second version of the new scheme above as DBDP1 and DBDP2, where the acronym DBDP stands for deep learning backward dynamic programming.

The intuition behind DBDP1 and DBDP2 is the following. For simplicity, take $f$ $=$ $0$, so that $F(t,x,y,z,h,\Delta)$ $=$ $y+z\trans\Delta$. The solution $u$ to the PDE \eqref{eq:PDEInit} should then approximately satisfy (see \eqref{uF})
\begin{align*}
u(t_{i+1}, X_{t_{i+1}}) & \approx  \; u(t_i,X_{t_i}) + D_x u(t_i,X_{t_i})\trans \sigma(t_i,X_{t_i}) \Delta W_{t_i}.  
\end{align*}
Consider the first scheme DBDP1, and suppose that at time $i+1$, $\widehat\Uc_{i+1}^{(1)}$ is an estimation of $u(t_{i+1,.})$. The quadratic  
loss function at time $i$ is then approximately equal to
\begin{align*}
\hat L_i^{(1)}(\theta) & \approx \; 
\E \Big| u(t_{i+1}, X_{t_{i+1}}) - \Uc_i(X_{t_i};\theta) - \Zc_i(X_{t_i};\theta)\trans\Delta W_{t_i} \Big|^2  \\
& \approx \; \E \Big[ \big| u(t_i,X_{t_i}) -  \Uc_i(X_{t_i};\theta) \big|^2   + 
\Delta t_i   \big|  \sigma\trans(t_i,X_{t_i}) D_x u(t_i,X_{t_i}) - \Zc_i(X_{t_i};\theta) \big|^2 \Big]. 
\end{align*}
Therefore, by minimizing over $\theta$ this quadratic loss function, via SGD based on simulations of $(X_{t_i},X_{t_{i+1}},\Delta W_{t_i})$ (called training data in the machine learning language), one expects the neural networks $\Uc_i$ and $\Zc_i$ to learn/approximate better and better the functions $u(t_i,.)$ and $\sigma\trans(t_i,)D_x u(t_i,)$ in view of the universal approximation theorem 
\cite{hor90universal}. Similarly, the second scheme 
DPDP2, which uses only neural network on the value functions, learns $u(t_i,.)$ by means of 
the neural network $\Uc_i$, and $\sigma\trans(t_i,)D_x u(t_i,)$ via  $\sigma\trans(t_i,)\hat D_x \Uc_i$. 
The rigorous arguments for the convergence of these schemes will be derived in the next section.  
 
\vspace{2mm}

The advantages of our two schemes, compared to the Deep BSDE algorithm, are the following:
\begin{itemize}
    \item by decomposing the global problem into smaller  ones, we may expect to help the gradient descent method to provide estimations closer to the real solution.
    The memory needed in \cite{han2017overcoming} can be a problem when taking too many time steps.
    \item at each time step, we initialize the weights and bias of the neural network to the weights and bias of the previous time step treated : this trick is commonly used in  iterative solvers of PDE, and allows us to start with a value close to the solution, hence avoiding local minima which are  too far away from the true solution. 
%methodology always used in iterative solvers in PDE methods permits to have a starting point close to the solution, and then to avoid local minima too far away from the true solution. 
 Besides the number of gradient iterations to achieve is rather small after the first resolution step. 
\end{itemize}
The small disadvantage is due to the Tensorflow structure. As it is done in python, the global graph creation takes much time  as it is repeated for each time step and the global resolution is a little bit time consuming : as the dimension of the problem increases, the time difference decreases and it becomes hard to compare the computational time for a given accuracy when the dimension is above 5.

\subsection{Extension to variational inequalities: scheme RDBDP}
\label{sec:varIneq}

Let us consider a variational inequality in the form
\begin{equation} \label{eq:IQV}
\left\{
\begin{aligned}
\min \big[ - \partial_t u - \Lc u - f(t,x,u,\sigma\trans D_x u) , u - g \big] & = 0 , \;\;\;\;\;\;\; t \in [0,T), \; x \in \R^d, \\
u(T,x) &=g(x), \;\;\; x\in\R^d.
\end{aligned}
\right.
\end{equation}
which arises, e.g., in optimal stopping problem and American option pricing in finance. It is known, see e.g. \cite{elk97}, that such variational inequality is related to reflected BSDE 
of the form
\begin{align} \label{RBSDE}
Y_t   &= \; g(\Xc_T) + \int_t^T f(s,\Xc_s,Y_s,Z_s) \diff s - \int_t^T Z_s\trans  \diff W_s + K_T - K_t,  \\
Y_t & \geq  \; g(X_t),  \;\;\; 0 \leq t \leq T,
\end{align}
where $K$ is an adapted non-decreasing process satisfying
\begin{align}
\int_0^T \big(Y_t - g(X_t) \big)  dK_t & = \; 0. 
\end{align}

The extension of our DBDP1 scheme for such variational inequality, and refereed to as RDBDP  scheme,  becomes
\begin{itemize}
\item  Initialize $\widehat\Uc_N$ $=$ $g$
\item  For $i$ $=$ $N-1,\ldots,0$, given $\widehat\Uc_{i+1}$, use a pair  of (multilayer) neural network $(\Uc_i(.;\theta),\Zc_i(.;\theta))$ $\in$ $\Nc\Nc_{d,1,L,m}^\varrho(\R^{N_m})\times\Nc\Nc_{d,d,L,m}^\varrho(\R^{N_m})$, and compute (by SGD) the minimizer of the expected quadratic loss function
\begin{equation} \label{eq:schemeVI}
\left\{ 
\begin{aligned} 
\hat L_i(\theta) & := \E \big| \widehat\Uc_{i+1}(X_{t_{i+1}}) - F(t_i,X_{t_i},\Uc_i(X_{t_i};\theta),\Zc_i(X_{t_i};\theta),\Delta t_i,\Delta W_{t_i}) \big|^2 \\
\theta_i^* & \in {\rm arg}\min_{\theta\in\R^{N_m}} \hat L_i(\theta). 
\end{aligned}
\right.
\end{equation}
Then, update: $\widehat\Uc_i$ $=$ $\max\big[\Uc_i(.;\theta_i^*),g]$, and set $\hat\Zc_i$ $=$ $\Zc(.;\theta_i ^*)$.  
\end{itemize}

%\begin{Remark}
%{\rm As for the scheme DBDP1, two neural networks can be used respectively for the approximation of $u(t_i,.)$ and for the approximation of $\sigma\trans(t_i,.)D_x u(t_i,.)$.
%}
%\ep
%\end{Remark}

\section{Convergence analysis}
\label{secconv} 

The main goal of this section is to prove  convergence of the DBDP schemes  towards the solution $(Y,Z)$ to the  BSDE \eqref{eqBSDE} (or reflected BSDE \eqref{RBSDE} for variational inequalities), and to provide a rate of convergence  that depends on the approximation errors  by neural networks.

\subsection{Convergence of DBDP1}

We assume the standard  Lipschitz conditions on $\mu$ and $\sigma$, which ensures the existence and uniqueness of an adapted solution $\Xc$ to the forward SDE \eqref{eq:SDE} satisfying for any $p$ $>$ $1$, 
\begin{align} \label{integX}
\E \big[ \sup_{0\leq t \leq T} |\Xc_t|^p \big] & \; <  C_p(1 + |x_0|^p),
\end{align} 
for some constant $C_p$ depending only on $p$, $b$, $\sigma$ and $T$. 
Moreover, we have the well-known error estimate with the Euler scheme $X$ $=$ $X^\pi$ defined in \eqref{eq:eulerSDE} with a time grid 
$\pi$ $=$ $\{t_0=0<t_1< \ldots < t_N = T\}$, with modulus $|\pi|$ s.t. $N|\pi|$ is bounded by a constant depending only on $T$ (hence independent of $N$): 
\begin{align} \label{estimEulerX}
\max_{i=0,\ldots,N-1} \E \Big[ |\Xc_{t_{i+1}} - X_{t_{i+1}}|^2 + \sup_{t\in[t_i,t_{i+1}]} | \Xc_t - X_{t_i}|^2 \Big] & = \; O(|\pi|). 
\end{align}
Here, the standard notation $O(|\pi|)$ means that 
$\limsup_{|\pi| \rightarrow 0}  \; |\pi|^{-1} O(|\pi|)$ $<$ $\infty$.

\vspace{1mm}

We shall make the standing usual assumptions on the driver  $f$ and the terminal data $g$.

\vspace{2mm}

\noindent \textbf{(H1)}  (i) There exists a constant $[f]_{_L}>0$ such that the driver $f$ satisfies:
\begin{equation}
    \left| f(t_2,x_2,y_2,z_2) - f(t_1,x_1,y_1,z_1)  \right| \leq [f]_{_L} \left( |t_2-t_1 |^{1/2} + |x_2-x_1| +|y_2-y_1| +|z_2-z_1| \right),
\end{equation}
for all $(t_1,x_1,y_1,z_1)$ and $(t_2,x_2,y_2,z_2)$ $\in [0,T] \times \R^d \times \R \times \R^d$.
Moreover, 
\begin{equation*}
\sup_{0 \leq t \leq T} |f(t,0,0,0)| < \infty.
\end{equation*}
(ii) The function $g$ satisfies a linear growth condition. 

\vspace{3mm}

%\noindent \textbf{(GA)}
%\marginpar{A t-on vraiment besoin de $g$ Lipschitz? En tout cas pas dans le Theoreme 3.2 et 3.3, Est ce dans Proposition 3.1?}
%The terminal condition $g$ is a Lipschitz function, i.e., there exists a constant $[g]_{_L}>0$ such that
%\begin{equation}
 %   \left| g(x_2)-g(x_1) \right| \leq [g]_{_L} \left| x_2-x_1 \right|, \;\;\; \forall x_1, x_2 \in \R^d.
%\end{equation}

%\vspace{2mm}

Recall that Assumption \textbf{(H1)} ensures the existence and uniqueness of an adapted solution $(Y,Z)$ to \eqref{eqBSDE} satisfying
\begin{align} \label{estiYZ}
\E \Big[ \sup_{0\leq t \leq T} |Y_t|^2 + \int_0^T |Z_t|^2 \diff t \Big] & < \; \infty. 
\end{align}
From the linear growth condition on $f$ in {\bf (H1)}, and \eqref{integX}, we also see that 
\begin{align} \label{integf2}
\E \Big[ \int_0^T |f(t,\Xc_t,Y_t,Z_t)|^2 \diff t \Big] &  < \; \infty. 
\end{align}
Moreover,  we  have the standard $L^2$-regularity result on $Y$:
\begin{align} \label{regulY}
\max_{i=0,\ldots,N-1} \E \Big[  \sup_{t\in[t_i,t_{i+1}]} | Y_t - Y_{t_i}|^2 \Big] & = \; O(|\pi|). 
\end{align} 
Let us also introduce the $L^2$-regularity of $Z$:
\begin{align}
\eps^Z(\pi) & := \; \E \bigg[ \sum_{i=0}^{N-1} \int_{t_i}^{t_{i+1}} |Z_t - \bar Z_{t_i}|^2 dt \bigg],
 \;\;\; \mbox{ with } \; \bar Z_{t_i} \; := \; \frac{1}{\Delta t_i} \E_i \Big[  \int_{t_i}^{t_{i+1}} Z_t dt \Big], 
\end{align}
where $\E_i$ denotes the conditional expectation given $\Fc_{t_i}$. Since $\bar Z$ is a  $L^2$-projection of $Z$, we know that $\eps^Z(\pi)$ converges to zero when $|\pi|$ goes to zero. Moreover, as shown in \cite{zhang04numerical}, when the terminal condition $g$ is also Lipschitz, we have 
\begin{align}
\eps^Z(\pi)  & = \; O(|\pi|). 
\end{align}

Let us  first investigate  the convergence of the scheme DBDP1 in \eqref{eq:scheme1}, and  define  (implicitly) 
\begin{equation} \label{defVCZ}
    \left\{ 
    \begin{array}{rcl}
\widehat\Vc_{t_i} & :=   & \E_i \big[ \widehat\Uc_{i+1}^{(1)}(X_{t_{i+1}}) \big] 
+ f(t_i,X_{t_i},\widehat\Vc_{t_i},\overline{{\widehat Z_{t_i}}}) \Delta t_i \\
\overline{{\widehat Z_{t_i}}} & := &  \frac{1}{\Delta t_i} \E_i\left[ \widehat\Uc_{i+1}^{(1)}(X_{t_{i+1}}) \Delta W_{t_i} \right],
\end{array}
\right. 
\end{equation}
for $i$ $=$ $0,\ldots,N-1$.   Notice that $\widehat\Vc_{t_i}$ is well-defined for $|\pi|$ small enough (recall that 
$f$ is Lipschitz) by a fixed point argument.  By the Markov property of the discretized forward process $(X_{t_i})_{i=0,\ldots,N}$, we note that there exists  
some deterministic functions $\hat v_i$ and $\overline{{\hat z_i}}$ s.t. 
\begin{align} \label{defhatv1} 
\widehat\Vc_{t_i}^{}  \; = \; \hat v_i^{}(X_{t_i}), & \mbox{ and } \;\;  
\overline{{\widehat Z_{t_i}}^{}} \; = \; \overline{{\hat z_i}^{}}(X_{t_i}), \;\;\;\;\;  i =0,\ldots,N-1.  
\end{align}
Moreover, by the martingale representation theorem, there exists an $\R^d$-valued square integrable process $(\widehat Z_t)_t$ such that 
\begin{align} \label{FBSDE}
\widehat\Uc_{i+1}^{(1)}(X_{t_{i+1}}) & = \;  \widehat\Vc_{t_i} -  
f(t_i,X_{t_i},\widehat\Vc_{t_i} ,\overline{{\widehat Z_{t_i}}}) \Delta t_i 
+ \int_{t_i}^{t_{i+1}}  \widehat Z_s\trans \diff W_s,  
%\\& = \;  F(t_i,X_{t_i},\widehat\Vc_{t_i},\overline{{\hat Z_{t_i}}},\Delta t_i,\Delta W_{t_i}) +  \int_{t_i}^{t_{i+1}} (\hat Z_s - \overline{{\hat Z_{t_i}}})\trans \diff W_s,  
%\;\;\;  i=0,\ldots,N-1. 
\end{align}
and by It\^o isometry, we have
\begin{align} \label{ZbarZ}
\overline{{\widehat Z_{t_i}}} &= \; \frac{1}{\Delta t_i} 
\E_i \Big[  \int_{t_i}^{t_{i+1}}  \widehat Z_s \diff s \Big], \;\;\;\;\;   i=0,\ldots,N-1. 
\end{align}

\vspace{2mm}

Let us now define a measure of the (squared) error for the DBDP1 scheme by
\begin{align}
\Ec\big[(\widehat\Uc^{(1)},\widehat\Zc^{(1)}),(Y,Z)\big] & := \; \max_{i=0,\ldots,N-1} \E \big|Y_{t_i}- \widehat\Uc_i^{(1)}(X_{t_i})\big|^2 + \E \bigg[ \sum_{i=0}^{N-1} \int_{t_i}^{t_{i+1}} 
\big| Z_t - \widehat\Zc_i^{(1)}(X_{t_i}) \big|^2 dt \bigg]. 
\end{align}

Our first main result  gives an error estimate of the  DBDP1 scheme   
in terms of the $L^2$-approximation errors of $\hat v_i$ and $\overline{{\hat z_i}}$ by neural networks $\Uc_i$ and $\Zc_i$, $i=0,\ldots,N-1$, assumed to be independent (see Remark \ref{remNN}), and defined as
\begin{align}
\eps_i^{\Nc,v} \; := \; \inf_{\xi} \E \big|\hat v_i(X_{t_i}) - \Uc_i(X_{t_i};\xi) \big|^2 ,
\hspace{7mm}  
 \eps_i^{\Nc,z} \; := \; \inf_{\eta} \E \big|\overline{{\hat z_i}^{}}(X_{t_i}) - 
\Zc_i(X_{t_i};\eta) \big|^2. 
\end{align}
Here, we fix the structure of the neural networks with input dimension $d$,  output dimension $d_1$ $=$ $1$ for $\Uc_i$, and $d_1$ $=$ $d$ for $\Zc_i$, number of layers $L$, and 
$m$ neurons for the hidden layers, and the parameters vary in the whole set $\R^{N_m}$ where $N_m$ is the number of parameters. 
From the universal approximation theorem (I) (\cite{horetal89}), we know that $\eps_i^{N N,v}$ and $\eps_i^{NN,z}$ converge to zero as $m$ goes to infinity, hence can be made 
arbitrary small for sufficiently large number of neurons.

\begin{Theorem} \emph{(Consistency of DBDP1)}
\label{theo:scheme1_1} Under {\bf (H1)},  
there exists a constant $C>0$, independent of $\pi$, such that
\begin{align}
\Ec\big[(\widehat\Uc^{(1)},\widehat\Zc^{(1)}),(Y,Z)\big] & \leq \; 
C  \Big(  \E \big|g(\Xc_{T}) - g(X_T) \big|^2 + |\pi| + \eps^Z(\pi) 
\\
& \hspace{9mm} 
+ \;  \sum_{i=0}^{N-1} \big(N \eps_i^{\Nc,v}  +   \eps_i^{\Nc,z}\big) \Big). 
\label{eq:theo1_scheme1}
\end{align}
\end{Theorem}

\begin{Remark}
{\rm The error contributions for the DBDP1 scheme in the r.h.s. of  estimation \eqref{eq:theo1_scheme1}  consists of four terms.  The first three terms correspond to the time discretization of BSDE, similarly as in \cite{bouchard2004discrete}, \cite{gobet2005regression},  
 namely (i) the strong approximation of the terminal condition (depending on the forward scheme and the terminal data $g$), and converging to zero, as $|\pi|$ goes to zero,  
with a rate $|\pi|$ when $g$ is Lipschitz by \eqref{estimEulerX} (see \cite{avi09} for irregular $g$), (ii) the strong approximation of the forward Euler scheme, and the $L^2$-regularity of $Y$, which gives a convergence of order $|\pi|$, (iii) the $L^2$-regularity of $Z$, which converges to zero, as $|\pi|$ goes to zero,  
with a rate $|\pi|$ when $g$ is Lipschitz. 
Finally, the better the neural networks are able to 
approximate/learn the functions $\hat v_i$ and $\overline{{\hat z_i}}$ at each time $i$ $=$ $0,\ldots,N-1$, 
the smaller is the last term in the error estimation. Moreover, given a prescribed accuracy for the neural network approximation error, the number of parameters of the employed deep neural networks grows at most polynomially in the PDE dimension, as recently proved in \cite{hutetal19} in the case of semi-linear heat equations.  
}
\ep
\end{Remark}

%\begin{Remark}
%{\rm
%The universal approximation theorem 
%\cite{hor90universal} is a priori  only valid on  compacts, but in  our case, some localization  methodology can be used as in \cite{chan2018machine} to show that  the error induced by the %network can be controlled at will.
%}
%\ep
%\end{Remark}

\vspace{3mm}

\noindent {\bf Proof of Theorem \ref{theo:scheme1_1}.}

\noindent In the following, $C$ will denote a positive generic constant independent of 
$\pi$, and that may take different values from line to line. 

\vspace{1mm}

\noindent {\it Step 1}. Fix $i$ $\in$ $\{0,\ldots,N-1\}$, and observe by \eqref{eqBSDE}, \eqref{defVCZ} that 
\begin{align} \label{reldifY}
\hspace{-5mm} Y_{t_i} - \widehat\Vc_{t_i} & =  \E_i\big[ Y_{t_{i+1}} - \widehat\Uc_{i+1}^{(1)}(X_{t_{i+1}}) \big] 
+ \E_i\Big[ \int_{t_i}^{t_{i+1}} f(t,\Xc_t,Y_t,Z_t) -  f(t_i,X_{t_i},\widehat\Vc_{t_i},\overline{{\widehat Z_{t_i}}}) \diff t \Big].  
\end{align}
By using Young inequality: $(a+b)^2$ $\leq$ $(1+\gamma \Delta t_i)a^2$ $+$ $(1+\frac{1}{\gamma \Delta t_i})b^2$ for some $\gamma$ $>$ $0$ to be chosen later, 
Cauchy-Schwarz inequality, the Lipschitz condition on $f$ in {\bf (H1)},  and the estimation \eqref{estimEulerX} on the forward process, we then have 
\begin{align}
\E\big| Y_{t_i} - \widehat\Vc_{t_i} \big|^2 & \leq \; 
(1 +\gamma\Delta t_i) \E \Big|  \E_i\big[ Y_{t_{i+1}} - \widehat\Uc_{i+1}^{(1)}(X_{t_{i+1}}) \big] \Big|^2 \\
& \;\;\; + 4 [f]^2_{_L} \Delta t_i \big(1+\frac{1}{\gamma \Delta t_i}\big) 
\Big\{ |\Delta t_i|^2 + \E\Big[ \int_{t_i}^{t_{i+1}} \big| Y_t - \widehat\Vc_{t_i} \big|^2 \diff t \Big] \\
& \hspace{4.3cm} + \; \E\Big[ \int_{t_i}^{t_{i+1}} \big| Z_t - \overline{{\widehat Z_{t_i}}} \big|^2 
\diff t \Big]   \Big\} \\
& \leq \;  (1 +\gamma\Delta t_i) \E \Big|  \E_i\big[ Y_{t_{i+1}} - \hat\Uc_{i+1}^{(1)}(X_{t_{i+1}}) \big] \Big|^2 \label{inegYVi} \\
& \;\;\; +  4 \frac{[f]^2_{_L}}{\gamma}  (1+ \gamma \Delta t_i)
\Big\{ C |\pi|^2 + 2 \Delta t_i \E\big| Y_{t_i} - \widehat\Vc_{t_i} \big|^2 
+ \E\Big[ \int_{t_i}^{t_{i+1}} \big| Z_t - \overline{{\widehat Z_{t_i}}} \big|^2  \diff t  \Big] \Big\},  
%\\
%& \hspace{3cm} + \;  
%\E \Big[\int_{t_i}^{t_{i+1}} \big| Z_t - \bar Z_{t_i}\big|^2 \diff t \Big] 
%+ \Delta t_i \E \big|\bar Z_{t_i} - \overline{{\widehat Z_{t_i}}} \big|^2  \Big] \Big\}, 
%\label{inegYVi}
\end{align} 
where we use in the last inequality the $L^2$-regularity \eqref{regulY} of $Y$. 
%Here the generic constant $C$ does not depend on $\gamma$. 

Recalling the definition of $\bar Z$ as a $L^2$-projection of $Z$, we observe that 
\begin{align} \label{Pythagore}
\E\Big[ \int_{t_i}^{t_{i+1}} \big| Z_t - \overline{{\widehat Z_{t_i}}} \big|^2  \diff t  \Big]  & = \; \E \Big[\int_{t_i}^{t_{i+1}} \big| Z_t - \bar Z_{t_i}\big|^2 \diff t \Big]  + 
\Delta t_i \E \big|\bar Z_{t_i} - \overline{{\widehat Z_{t_i}}} \big|^2.  
\end{align}
By multiplying  equation \eqref{eqBSDE}  between $t_i$ and $t_{i+1}$ by $\Delta W_{t_i}$, and  using It\^o isometry, we have together with \eqref{defVCZ} 
\begin{align}
\Delta t_i \big( \bar Z_{t_i} - \overline{{\widehat Z_{t_i}}} \big) & = 
\; \E_i \big[ \Delta W_{t_i}  \big( Y_{t_{i+1}} - \widehat\Uc_{i+1}^{(1)}(X_{t_{i+1}}) \big) \big]    + \E_i \Big[ \Delta W_{t_i}  \int_{t_i}^{t_{i+1}} f(t,\Xc_t,Y_t,Z_t) \diff t  \Big] \\
& = \;  \E_i \Big[ \Delta W_{t_i}  \Big( Y_{t_{i+1}} - \widehat\Uc_{i+1}^{(1)}(X_{t_{i+1}})  -  \E_i\big[ Y_{t_{i+1}} - \widehat\Uc_{i+1}^{(1)}(X_{t_{i+1}})\big]    \Big) \Big]   \\
& \;\;\;\; + \; \E_i \Big[ \Delta W_{t_i}  \int_{t_i}^{t_{i+1}} f(t,\Xc_t,Y_t,Z_t) \diff t  \Big]. 
\end{align}
By Cauchy-Schwarz inequality, and law of iterated conditional expectations,  this implies
\begin{align}
\Delta t_i \E \big| \bar Z_{t_i} - \overline{{\widehat Z_{t_i}}} \big|^2 & \leq \;  2 d \Big( \E \big|Y_{t_{i+1}} - \widehat\Uc_{i+1}^{(1)}(X_{t_{i+1}}) \big|^2 - 
\E \Big| \E_i\big[ Y_{t_{i+1}} - \widehat\Uc_{i+1}^{(1)}(X_{t_{i+1}})\big]   \Big|^2  \Big) \\
& \;\;\;\;\; + 2 d \Delta t_i \E \Big[ \int_{t_i}^{t_{i+1}} |f(t,\Xc_t,Y_t,Z_t)|^2  \diff t  \Big]. \label{inegZi} 
\end{align}
Then,  by  plugging \eqref{Pythagore} and \eqref{inegZi}  into \eqref{inegYVi}, and choosing $\gamma$ $=$ $8 d [f]^2_{_L}$, we have
\begin{align}
\E\big| Y_{t_i} - \widehat\Vc_{t_i} \big|^2 & \leq \;  C \Delta t_i \E\big| Y_{t_i} - \widehat\Vc_{t_i} \big|^2  
+ (1 +\gamma\Delta t_i) \E \big|Y_{t_{i+1}} - \widehat\Uc_{i+1}^{(1)}(X_{t_{i+1}}) \big|^2 + C |\pi|^2 \\
& \;\;\; + \; C \E \Big[\int_{t_i}^{t_{i+1}} \big| Z_t - \bar Z_{t_i}\big|^2 \diff t \Big]  + C \Delta t_i \E \Big[ \int_{t_i}^{t_{i+1}} |f(t,\Xc_t,Y_t,Z_t)|^2  \diff t  \Big],
\end{align}
and thus for $|\pi|$ small enough: 
\begin{align}
\E\big| Y_{t_i} - \widehat\Vc_{t_i} \big|^2 & \leq \; (1 + C |\pi|) \E \big|Y_{t_{i+1}} - \widehat\Uc_{i+1}^{(1)}(X_{t_{i+1}}) \big|^2 +  C |\pi|^2  \\
& \hspace{.2cm}+  C \E \Big[\int_{t_i}^{t_{i+1}} \big| Z_t - \bar Z_{t_i}\big|^2 \diff t \Big]  + C |\pi| \E \Big[ \int_{t_i}^{t_{i+1}} |f(t,\Xc_t,Y_t,Z_t)|^2  \diff t  \Big]. 
\label{interYUV}
\end{align}

\vspace{1mm}

\noindent {\it Step 2.} By using Young inequality in the form: $(a+b)^2$ $\geq$ $(1-  |\pi|)a^2$ $+$ $(1-\frac{1}{|\pi|})b^2$ $\geq$ 
$(1-  |\pi|)a^2$ $-$ $\frac{1}{ |\pi|}b^2$, we have 
\begin{align}
\E\big| Y_{t_i} - \widehat\Vc_{t_i} \big|^2 & = \E\big| Y_{t_i} -  \widehat\Uc_{i}^{(1)}(X_{t_{i}}) + \widehat\Uc_{i}^{(1)}(X_{t_{i}}) - \widehat\Vc_{t_i} \big|^2 \\
& \geq \;  (1-  |\pi|) \E\big| Y_{t_i} -  \widehat\Uc_{i}^{(1)}(X_{t_{i}}) \big|^2 - \frac{1}{ |\pi|} \E \big| \widehat\Uc_{i}^{(1)}(X_{t_{i}}) - \widehat\Vc_{t_i} \big|^2.  \label{YUinter}
\end{align}
By plugging this last inequality into \eqref{interYUV}, we then get for $|\pi|$ small enough 
\begin{align}
\E\big| Y_{t_i} - \widehat\Uc_{i}^{(1)}(X_{t_{i}}) \big|^2 & \leq \; (1 + C |\pi|) \E \big|Y_{t_{i+1}} - \widehat\Uc_{i+1}^{(1)}(X_{t_{i+1}}) \big|^2 +  C |\pi|^2  \\
& \;\;\;\;\; + \; C \E \Big[\int_{t_i}^{t_{i+1}} \big| Z_t - \bar Z_{t_i}\big|^2 \diff t \Big]  + C |\pi| \E \Big[ \int_{t_i}^{t_{i+1}} |f(t,\Xc_t,Y_t,Z_t)|^2  \diff t  \Big] \\
& \;\;\;\;\; + \;  C N  \E \big| \widehat\Vc_{t_i}  - \widehat\Uc_{i}^{(1)}(X_{t_{i}}) \big|^2. 
\end{align}
From discrete Gronwall's lemma (or by induction), and recalling the terminal condition $Y_{t_N}$ $=$ $g(\Xc_T)$, $\widehat\Uc_{i}^{(1)}(X_{t_{N}})$ $=$ $g(X_T)$, 
the definition $\eps^Z(\pi)$  of the $L^2$-regularity of $Z$, and \eqref{integf2}, this yields
\begin{align}
\max_{i=0,\ldots,N-1} \E\big| Y_{t_i} - \widehat\Uc_{i}^{(1)}(X_{t_{i}}) \big|^2 & \leq \; C  \E \big|g(\Xc_{T}) - g(X_T) \big|^2 + C |\pi| + C \eps^Z(\pi)  \\
& \hspace{5mm}  + \;  C N \sum_{i=0}^{N-1}  \E \big| \widehat\Vc_{t_i}  -  \widehat\Uc_{i}^{(1)}(X_{t_{i}}) \big|^2. \label{estimYinter} 
\end{align}

\vspace{1mm}

\noindent {\it Step 3.} Fix  $i$ $\in$ $\{0,\ldots,N-1\}$. By using relation \eqref{FBSDE} in the expression of the expected quadratic loss function in \eqref{eq:scheme1}, and recalling the definition of $\overline{{\widehat Z_{t_i}}}$ as a $L^2$-projection of $\widehat Z_t$, we have 
for all parameters $\theta$ $=$ $(\xi,\eta)$ of the neural networks $\Uc_i(.;\xi)$ and $\Zc_i(.;\eta)$
\begin{align} \label{LtildeL}
\hat L_i^{(1)}(\theta) &= \;  \tilde L_i(\theta) + \E \Big[ \int_{t_i}^{t_{i+1}} \big| \widehat Z_t -  \overline{{\widehat Z_{t_i}}} \big|^2 \diff t \Big]
\end{align}
with 
\begin{align}
\tilde L_i(\theta) & := \; \E \Big| \widehat\Vc_{t_i}  - \Uc_i(X_{t_i};\xi) + 
\big( f(t_i,X_{t_i},\Uc_i(X_{t_i};\xi),\Zc_i(X_{t_i};\eta)) - f(t_i,X_{t_i},\widehat\Vc_{t_i},\overline{{\widehat Z_{t_i}}}) \big) \Delta t_i \Big|^2 \\
& \hspace{.6cm}  + \;  \Delta t_i \E \big| \overline{{\widehat Z_{t_i}}} - \Zc_i(X_{t_i};\eta) \big|^2.
\end{align}
By using Young inequality: $(a+b)^2$ $\leq$ $(1+\gamma \Delta t_i)a^2$ $+$ $(1+\frac{1}{\gamma \Delta t_i})b^2$, together with the Lipschitz condition on $f$ in {\bf (H1)}, we clearly see that 
\begin{align} \label{L<}
\tilde L_i(\theta) & \leq \; (1 + C \Delta t_i) \E \big| \widehat\Vc_{t_i}  - \Uc_i(X_{t_i};\xi) \big|^2 + C \Delta t_i \E \big| \overline{{\widehat Z_{t_i}}} - \Zc_i(X_{t_i};\eta) \big|^2.  
\end{align} 
On the other hand, using Young inequality in the form: $(a+b)^2$ $\geq$ $(1- \gamma \Delta t_i)a^2$ $+$ $(1-\frac{1}{\gamma \Delta t_i})b^2$ $\geq$ 
$(1- \gamma \Delta t_i)a^2$ $-$ $\frac{1}{\gamma \Delta t_i}b^2$, together with the Lipschitz condition on $f$, we have
\begin{align}
\tilde L_i(\theta) & \geq \; (1 - \gamma \Delta t_i) \E \big| \widehat\Vc_{t_i}  - \Uc_i(X_{t_i};\xi) \big|^2 
- \frac{2 \Delta t_i [f]^2_{_L}}{\gamma} \Big(  \E \big| \widehat\Vc_{t_i}  - \Uc_i(X_{t_i};\xi) \big|^2  + \E \big| \overline{{\widehat Z_{t_i}}} - \Zc_i(X_{t_i};\eta) \big|^2 \Big) \nonumber \\
& \hspace{.6cm} + \;   \Delta t_i \E \big| \overline{{\widehat Z_{t_i}}} - \Zc_i(X_{t_i};\eta) \big|^2.  
\end{align} 
By choosing $\gamma$ $=$ $4[f]_{_L}^2$, this yields
\begin{align} \label{L>}
\tilde L_i(\theta) & \geq \; (1 - C \Delta t_i) \E \big| \widehat\Vc_{t_i}  - \Uc_i(X_{t_i};\xi) \big|^2 +   \frac{\Delta t_i}{2} \E \big| \overline{{\widehat Z_{t_i}}} - \Zc_i(X_{t_i};\eta) \big|^2. 
\end{align} 

\vspace{1mm}

\noindent {\it Step 4.}  Fix  $i$ $\in$ $\{0,\ldots,N-1\}$, and 
take $\theta_i^*$ $=$ $(\xi_i^*,\eta_i^*)$ $\in$ ${\rm arg}\min_\theta \hat L_i^{(1)}(\theta)$ so that $\widehat\Uc_i^{(1)}$ $=$ $\Uc_i(.;\xi_i^*)$, and 
$\widehat\Zc_i^{(1)}$ $=$ $\Zc_i(.;\eta_i^*)$. By \eqref{LtildeL}, notice that $\theta_i^*$ $\in$ ${\rm arg}\min_\theta \tilde L_i(\theta)$. From \eqref{L>} and \eqref{L<}, we then have for all 
$\theta$ $=$ $(\xi,\eta)$
\begin{align}
(1 - C \Delta t_i) \E \big| \widehat\Vc_{t_i}  - \widehat\Uc_i^{(1)}(X_{t_i}) \big|^2 +   \frac{\Delta t_i}{2} \E \big| \overline{{\widehat Z_{t_i}}} - \widehat\Zc_i^{(1)}(X_{t_i}) \big|^2 \\
& \hspace{-8cm} \leq \tilde L_i(\theta_i^*) \; \leq \tilde L_i(\theta)  \; \leq \; (1 + C \Delta t_i) \E \big| \widehat\Vc_{t_i}  - \Uc_i(X_{t_i};\xi) \big|^2 + C \Delta t_i \E \big| \overline{{\widehat Z_{t_i}}} - \Zc_i(X_{t_i};\eta) \big|^2. \nonumber 
\end{align}
For $|\pi|$ small enough, and recalling \eqref{defhatv1}, this implies
\begin{align} \label{estimVU}
\E \big| \widehat\Vc_{t_i}  - \widehat\Uc_i^{(1)}(X_{t_i}) \big|^2 +   \Delta t_i \E \big| \overline{{\widehat Z_{t_i}}} - \widehat\Zc_i^{(1)}(X_{t_i}) \big|^2
& \leq \; C \eps_i^{\Nc,v} + C \Delta t_i \eps_i^{\Nc,z}. 
\end{align}
Plugging this last inequality into \eqref{estimYinter}, we obtain
\begin{align} 
\max_{i=0,\ldots,N-1} \E\big| Y_{t_i} - \widehat\Uc_{i}^{(1)}(X_{t_{i}}) \big|^2 & \leq \; C  \E \big|g(\Xc_{T}) - g(X_T) \big|^2 + C |\pi| + C \eps^Z(\pi)  \\
& \hspace{6mm}  + \;  C  \sum_{i=0}^{N-1}   \big( N \eps_i^{\Nc,v} +  \eps_i^{\Nc,z} \big), \label{estimYfin} 
\end{align} 
which proves the consistency of the $Y$-component in \eqref{eq:theo1_scheme1}.   
 
\vspace{1mm}

\noindent {\it Step 5.}  Let us finally  prove the consistency of the $Z$-component. From \eqref{Pythagore} and  \eqref{inegZi}, we have for any $i$ $=$ $0,\ldots,N-1$:
\begin{align}
\E\Big[ \int_{t_i}^{t_{i+1}} \big| Z_t - \overline{{\widehat Z_{t_i}}}  \big|^2  \diff t  \Big]  & \leq \; 
\E\Big[ \int_{t_i}^{t_{i+1}} \big| Z_t - \bar Z_{t_i} \big|^2  \diff t  \Big] +  2d  |\pi| \E \Big[ \int_{t_i}^{t_{i+1}} |f(t,\Xc_t,Y_t,Z_t)|^2  \diff t  \Big] \\
& \hspace{.6cm} + \; 2d \Big( \E \big|Y_{t_{i+1}} - \widehat\Uc_{i+1}^{(1)}(X_{t_{i+1}}) \big|^2 - 
\E \Big| \E_i\big[ Y_{t_{i+1}} - \widehat\Uc_{i+1}^{(1)}(X_{t_{i+1}})\big]   \Big|^2  \Big) 
\end{align}
By summing over $i$ $=$ $0,\ldots,N-1$, we get (recall \eqref{integf2}) 
\begin{align}
\E \Big[ \sum_{i=0}^{N-1} \int_{t_i}^{t_{i+1}} \big| Z_t - \overline{{\widehat Z_{t_i}}}  \big|^2  \diff t  \Big] & \leq \;  \eps^Z(\pi) + C |\pi| + 2d \E \big|g(\Xc_{T}) - g(X_T) \big|^2 \label{Zinter} \\
& \; +  \: 2d \sum_{i=0}^{N-1} \Big( \E \big|Y_{t_{i}} - \widehat\Uc_{i}^{(1)}(X_{t_{i}}) \big|^2 - 
\E \Big| \E_i\big[ Y_{t_{i+1}} - \widehat\Uc_{i+1}^{(1)}(X_{t_{i+1}})\big]   \Big|^2 \Big) \nonumber
\end{align}
where we change the indices in the last summation. Now, from \eqref{inegYVi}, \eqref{YUinter}, we have
 \begin{align}
2d \Big( \E \big|Y_{t_{i}} - \widehat\Uc_{i}^{(1)}(X_{t_{i}}) \big|^2 -  \E \Big| \E_i\big[ Y_{t_{i+1}} - \widehat\Uc_{i+1}^{(1)}(X_{t_{i+1}})\big]   \Big|^2  \Big) \\
& \hspace{-7cm}\leq  
\Big(\frac{1 +\gamma |\pi|}{1  - |\pi|}  - 1 \Big)\E \Big|  \E_i\big[ Y_{t_{i+1}} - \hat\Uc_{i+1}^{(1)}(X_{t_{i+1}}) \big] \Big|^2  \\
&  \hspace{-6.4cm}+   \;  \frac{8d[f]^2_{_L}}{\gamma}  \frac{1+ \gamma |\pi|}{1 -  |\pi|} 
\Big\{ C |\pi|^2  +  |\pi|  \E\big| Y_{t_i} - \widehat\Vc_{t_i} \big|^2 
+ \E\Big[ \int_{t_i}^{t_{i+1}} \big| Z_t - \overline{{\widehat Z_{t_i}}} \big|^2  \diff t  \Big] \Big\}  \\
& \hspace{-6.4cm}+ \; \frac{2d}{|\pi|(1-|\pi|)} \E \big| \widehat\Uc_{i}^{(1)}(X_{t_{i}}) - \widehat\Vc_{t_i} \big|^2. 
 \end{align}
 We now choose $\gamma$ $=$ $24d[f]^2_{_L}$ so that $\frac{8d[f]^2_{_L}}{\gamma}  (1+ \gamma |\pi|)/(1 -  |\pi|)$ $\leq$ $1/2$ for $|\pi|$ small enough, and by plugging into \eqref{Zinter}, 
 we obtain (note also that $\big[(1 +\gamma |\pi|)/(1 -  |\pi|)  - 1 \big]$ $=$ $O(|\pi|)$):
 \begin{align}
\frac{1}{2} \E \Big[ \sum_{i=0}^{N-1} \int_{t_i}^{t_{i+1}} \big| Z_t - \overline{{\widehat Z_{t_i}}}  \big|^2  \diff t  \Big] & \leq \; 
\eps^Z(\pi) + C |\pi| +  C \max_{i=0,\ldots,N} \E\big| Y_{t_i} - \widehat\Uc_{i}^{(1)}(X_{t_{i}}) \big|^2 \\
& \hspace{.6cm}+ \; \frac{1}{2} |\pi| \sum_{i=0}^{N-1} \E\big| Y_{t_i} - \widehat\Vc_{t_i} \big|^2  + CN \sum_{i=0}^{N-1} \E \big| \widehat\Uc_{i}^{(1)}(X_{t_{i}}) - \widehat\Vc_{t_i} \big|^2 \\
& \leq \;  C \eps^Z(\pi) + C |\pi| +  C \max_{i=0,\ldots,N} \E\big| Y_{t_i} - \widehat\Uc_{i}^{(1)}(X_{t_{i}}) \big|^2 \\
& \hspace{.6cm} + \; CN \sum_{i=0}^{N-1} \E \big| \widehat\Uc_{i}^{(1)}(X_{t_{i}}) - \widehat\Vc_{t_i} \big|^2 \\
& \leq \; C  \E \big|g(\Xc_{T}) - g(X_T) \big|^2 + C |\pi| + C \eps^Z(\pi)  \\
& \hspace{6mm}  + \;  C  \sum_{i=0}^{N-1}   \big( N \eps_i^{\Nc,v} +  \eps_i^{\Nc,z} \big), \label{estimZhatZ}
\end{align}
 where we used \eqref{interYUV} and \eqref{integf2}  in the second inequality, and \eqref{estimVU} and \eqref{estimYfin} in the last inequality. 
 
 By writing that 
 \begin{align}
 \E\Big[ \int_{t_i}^{t_{i+1}} \big| Z_t - \widehat\Zc_i^{(1)}(X_{t_i}) \big|^2  \diff t  \Big] & \leq \;  2 \E\Big[ \int_{t_i}^{t_{i+1}} \big| Z_t - \overline{{\widehat Z_{t_i}}}  \big|^2  \diff t  \Big] 
 + 2 \Delta t_i  \E \big| \overline{{\widehat Z_{t_i}}} - \widehat\Zc_i^{(1)}(X_{t_i}) \big|^2,
 \end{align}
 and using  \eqref{estimVU}, \eqref{estimZhatZ}, we obtain after summation over $i$ $=$ $0,\ldots,N-1$,  the required error estimate for the $Z$-component as in \eqref{estimYfin}, and this 
ends the proof.    
\ep

\subsection{Convergence of DBDP2}
%\marginpar{XW total variation smaller than $\gamma_m$? "total variation" %prête à confusion. C'est ta def il me semble. Faire un renvoi à l'équation %dans section 2 ? }
 We shall consider neural networks with one hidden layer, $m$ neurons with total variation smaller than $\gamma_m$ (see Section \ref{secNN}),  a $C^3$ activation function $\varrho$ with linear growth condition, and bounded derivatives,  e.g.,  a sigmoid activation function, or a $\tanh$ function: this class of neural networks is then represented  by the parametric set of functions 
 \begin{align}
 \Nc\Nc_{d,1,2,m}^\varrho(\Theta_m^\gamma) 
 &:= \left\{ x \in \R^d \mapsto \mathcal{U}(x;\theta) = \sum_{i=1}^{m} c_i \varrho(a_i.x+b_i) + b_0, \; \theta=(a_i,b_i,c_i,b_0)_{i=1}^m \; \in \; \Theta^\gamma_m  \right\}, \nonumber
 \end{align}
 with
  \begin{align}
 \Theta^\gamma_m &:= \left\{ \theta=(a_i,b_i,c_i,b_0)_{i=1}^m:  \;  \max_{i=1,\ldots,m} |a_i| \leq \gamma_m,   \;  \sum_{i=1}^{m}  |c_i| \leq \gamma_m \right\},  
 \end{align}
for some sequence $(\gamma_m)_m$ converging to $\infty$, as $m$ goes to infinity, 
and such that  
\begin{equation} \label{convgamma}
\begin{array}{cc}
\frac{\gamma_{m}^6}{N}  \xrightarrow[m,N \to \infty]{} 0.
\end{array}
\end{equation}
%In the sequel, we denote by $\Theta^N$ the set of parameters %$\Theta^{\gamma_{_N}}_{m_{_N}}$, and by $\Nc\Nc^N$ the associated set of neural networks 
%$\Nc\Nc_{d,1,2,m_{_N}}^\varrho(\Theta^N)$.  
Notice that the neural networks in $\Nc\Nc_{d,1,2,m}^\varrho(\Theta_m^\gamma)$  have their first, second and third derivatives uniformly bounded w.r.t. the state variable $x$. 
More precisely, there exists some constant $C$ depending only on  $d$ and  the derivatives of $\varrho$ s.t. for any $\Uc$ $\in$ $\Nc\Nc_{d,1,2,m}^\varrho(\Theta_m^\gamma)$, 
\begin{equation} \label{bounderivU} 
\left\{
\begin{aligned}
\sup_{x\in\R^d,\theta\in\Theta_m^\gamma}  \Big| D_x \mathcal{U}(x;\theta) \Big|  \; \leq \; C   \gamma_{m}^2, \quad  
\sup_{x\in\R^d,\theta\in\Theta_m^\gamma}   \Big| D_x^2  \mathcal{U}(x;\theta) \Big|  \; \leq  \; C \gamma_{m}^3,   \\
\text{ and } \;\;\;  \sup_{x\in\R^d,\theta\in\Theta_m^\gamma} 
\Big| D_x^3  \mathcal{U}(x;\theta) \Big|  \; \leq  \; C \gamma_{m}^4. 
 \end{aligned}
 \right.
 \end{equation}

\vspace{1mm}

Let us  investigate  the convergence of the scheme DBDP2 in \eqref{eq:scheme2} with neural networks in $\Nc\Nc_{d,1,2,m}^\varrho(\Theta_m^\gamma)$, 
and  define  for $i$ $=$ $0,\ldots,N-1$:  
\begin{equation} \label{defVCZ2}
    \left\{ 
    \begin{array}{rcl}
\widehat\Vc_{t_i} & :=   & \E_i \big[ \widehat\Uc_{i+1}^{(2)}(X_{t_{i+1}}) \big] 
+ f(t_i,X_{t_i},\widehat\Vc_{t_i},\overline{{\widehat Z_{t_i}}}) \Delta t_i \; = \;  \hat v_i^{}(X_{t_i}),  \\
\overline{{\widehat Z_{t_i}}} & := &  \frac{1}{\Delta t_i} \E_i\left[ \widehat\Uc_{i+1}^{(2)}(X_{t_{i+1}}) \Delta W_{t_i} \right] \; = \; \overline{{\hat z_i}^{}}(X_{t_i}). 
\end{array}
\right. 
\end{equation}

\vspace{1mm}

A  measure of the (squared) error for the DBDP2 scheme is defined similarly as in DBDP1 scheme: 
\begin{align}
\Ec\big[(\widehat\Uc^{(2)},\widehat\Zc^{(2)}),(Y,Z)\big] & := \; 
\max_{i=0,\ldots,N-1} \E \big|Y_{t_i}- \widehat\Uc_i^{(2)}(X_{t_i})\big|^2 + \E \bigg[ \sum_{i=0}^{N-1} \int_{t_i}^{t_{i+1}} 
\big| Z_t - \widehat\Zc_i^{(2)}(X_{t_i}) \big|^2 dt \bigg]. 
\end{align}

%\vspace{1mm}

Our second main result  gives an error estimate of the  DBDP2 scheme   
in terms of the $L^2$-approximation errors of $\hat v_i$ and its derivative (which exists under assumption detailed below)  by neural networks $\Uc_i$ $\in$ $\Nc\Nc_{d,1,2,m}^\varrho(\Theta_m^\gamma)$, $i=0,\ldots,N-1$, and defined as
\begin{align}
\eps_i^{\Nc,m} \; := \; \inf_{\theta \in \Theta_m^\gamma}\Big\{  \E \big|\hat v_i(X_{t_i}) - \Uc_i(X_{t_i};\theta) \big|^2  
+ \Delta t_i \E \big| \sigma\trans(t_i,X_{t_i}) \big(  D_x \hat v_i(X_{t_i}) - D_x \Uc_i(X_{t_i};\theta) \big) \big|^2  \Big\}, \nonumber
\end{align}
which are expected to be small in view of the universal approximation theorem (II), see discussion in Remark \ref{remNN_approxError}.

\vspace{1mm}

We also require the additional  conditions on the coefficients:
 
\vspace{2mm}

 \noindent {\bf (H2)}  (i) The functions $x$ $\mapsto$ $\mu(t,.)$, $\sigma(t,.)$ are $C^1$ with bounded derivatives uniformly w.r.t. $(t,x)$ $\in$ $[0,T]\times\R^d$. 

\vspace{1mm}

\noindent (ii) The function  $(x,y,z)$ $\mapsto$ $f(t,.)$ is $C^1$  with bounded derivatives uniformly w.r.t. $(t,x,y,z)$  in $[0,T]\times \R^d\times\R\times\R^d$.

\vspace{2mm}

\begin{Theorem} \label{theoconv2} \emph{(Consistency of DBDP2)}
\label{theo:scheme2} Under {\bf (H1)}-{\bf (H2)},  there exists a constant $C>0$, independent of $\pi$, such that
\begin{align}
\Ec\big[(\widehat\Uc^{(2)},\widehat\Zc^{(2)}),(Y,Z)\big] & \leq \; 
C  \Big(  \E \big|g(\Xc_{T}) - g(X_T) \big|^2 + \frac{\gamma_{m}^6}{N}  + \eps^Z(\pi) + N \sum_{i=0}^{N-1} \eps_i^{\Nc,m}  \Big). 
\label{eq:theo2_scheme2}
\end{align}
\end{Theorem}
\noindent {\bf Proof.} For simplicity of notations, we assume $d$ $=$ $1$, and only detail the arguments that differ from the proof of Theorem \ref{eq:theo1_scheme1}. 
From \eqref{defVCZ2}, and the Euler scheme \eqref{eq:eulerSDE}, we have 
\begin{align}
\hat v_i(x) &  =  \; \tilde v_i(x) + \Delta t_i  f(t_i,x,\hat v_i(x), \overline{{\hat z_i}}(x)), \;\;\;  \tilde v_i(x) \; := \; \E\big[ \hat u_{i+1}(X_{t_{i+1}}^x) \big], \;\;  x \in \R^d,  \\
\overline{{\hat z_i}}(x) &  = \; \frac{1}{\Delta t_i} \E \big[ \hat u_{i+1}(X_{t_{i+1}}^x) \Delta W_{t_i} \big], \;\;\; X_{t_{i+1}}^x \; = \; x + \mu(t_i,x) \Delta t_i + \sigma(t_i,x) \Delta W_{t_i}. 
\end{align}
Under assumption {\bf (H2)}(i), and recalling that $\hat u_{i+1}$ $=$ $\Uc_{i+1}(.;\theta_{i+1}^*)$  is $C^2$ with bounded derivatives, we see that $\tilde v_i$ is $C^1$ with
\begin{align}
D_x \tilde v_i(x) & = \;   \E \Big[ \big(1 + D_x \mu(t_i,x) \Delta t_i  + D_x \sigma(t_i,x) \Delta W_{t_i} \big) D_x \hat u_{i+1}(X_{t_{i+1}}^x) \Big] \\
& = \; \E \big[ D_x \hat u_{i+1}(X_{t_{i+1}}^x) \big]  +  \Delta t_i  \; R_i(x)  \label{derivtildev} \\
R_i(x) & := \;  D_x \mu(t_i,x) \E \big[ D_x \hat u_{i+1}(X_{t_{i+1}}^x) \big]  + \sigma(t_i,x) D_x \sigma(t_i,x)  \E \big[ D_x^2 \hat u_{i+1}(X_{t_{i+1}}^x) \big],   
\end{align}
where we use integration by parts in the second equality. Similarly, we have
\begin{equation}
\left\{
\begin{aligned}
\overline{{\hat z_i}}(x) &  = \; \sigma(t_i,x)\E\big[ D_x \hat u_{i+1}(X_{t_{i+1}}^x) \big],  \label{IPPZ} \\
D_x \overline{{\hat z_i}}(x)  & = \;  D_x \sigma(t_i,x)\E\big[ D_x \hat u_{i+1}(X_{t_{i+1}}^x) \big]  + \sigma(t_i,x)\E\big[ D_x^2 \hat u_{i+1}(X_{t_{i+1}}^x) \big] + 
\Delta t_i  \; \sigma(t_i,x) G_i(x) \\
G_i(x) & := \; D_x \mu(t_i,x) \E \big[ D_x^2 \hat u_{i+1}(X_{t_{i+1}}^x) \big]  + \sigma(t_i,x) D_x \sigma(t_i,x)  \E \big[ D_x^3 \hat u_{i+1}(X_{t_{i+1}}^x) \big].   
\end{aligned}
\right.
\end{equation} 
Denoting by $\hat f_i(x)$ $=$ $f(t_i,x,\hat v_i(x), \overline{{\hat z_i}}(x))$, it follows by the implicit function theorem, and for $|\pi|$ small enough, that $\hat v_i$ is $C^1$ with derivative given by
\begin{align}
D_x \hat v_i(x) & = \;  D_x \tilde v_i(x)  + \Delta t_i \Big( D_x \hat f_i(x)  +  D_y \hat f_i(x) D_x \hat v_i(x) + D_z \hat f_i(x) D_x \overline{{\hat z_i}}(x)  \Big)
\end{align}
and thus by \eqref{derivtildev}-\eqref{IPPZ} 
\begin{align}
\big( 1 - \Delta t_i D_y \hat f_i(x) \big) \sigma(t_i,x) D_x \hat v_i(x)  & =  \;  \overline{{\hat z_i}}(x)   +  \Delta t_i \sigma(t_i,x) \Big( R_i(x) +  D_x \hat f_i(x) +  D_z \hat f_i(x) D_x \overline{{\hat z_i}}(x)  \Big).  
\end{align}
Under {\bf (H2)}, by the linear growth condition on $\sigma$, and using the bounds on the derivatives of the neural networks in $\Nc\Nc_{d,1,2,m}^\varrho(\Theta_m^\gamma)$  in \eqref{bounderivU}, we then have 
\begin{align} \label{interzv}
\E \Big|  \sigma(t_i,X_{t_i}) D_x \hat v_i(X_{t_i})  -  \overline{{\widehat Z_{t_i}}} \Big|^2 & \leq \; C( \gamma_{m}^6 + |\pi|^2 \gamma_{m}^8) |\pi|^2. 
\end{align}
Next, by  the same arguments as in Steps 3 and  4 in the proof of  Theorem \ref{theo:scheme1_1} (see in particular \eqref{estimVU}), we have for $|\pi|$ small enough, 
\begin{align} 
\E \big| \widehat\Vc_{t_i}  - \widehat\Uc_i^{(2)}(X_{t_i}) \big|^2 +   \Delta t_i \E \big| \overline{{\widehat Z_{t_i}}} - \widehat\Zc_i^{(2)}(X_{t_i}) \big|^2 \\
& \hspace{-5cm}\leq  C \E\big[ \big|\hat v_i(X_{t_i}) - \Uc_i(X_{t_i};\theta) \big|^2 \big]  + C \Delta t_i 
 \E  \big| \overline{{\widehat Z_{t_i}}}  -   \sigma(t_i,X_{t_i}) \hat D_x \Uc_i(X_{t_i};\theta) \big|^2, 
\end{align}
for all $\theta$ $\in$ $\Theta^N$, and then with \eqref{interzv}, and by definition of $\eps_i^{NN,v,2}$: 
\begin{align} \label{estimU2}
 \E \big| \widehat\Vc_{t_i}  - \widehat\Uc_i^{(2)}(X_{t_i}) \big|^2 +   \Delta t_i \E \big| \overline{{\widehat Z_{t_i}}} - \widehat\Zc_i^{(2)}(X_{t_i}) \big|^2 
 & \leq \;  C \eps_i^{NN,v,2} + C ( \gamma_{m}^6 + |\pi|^2 \gamma_{m}^8)   |\pi|^3.  
\end{align}
On the other hand, by the same arguments as in Steps 1 and  2 in the proof of  Theorem \ref{theo:scheme1_1} (see in particular \eqref{estimYinter}), we have 
\begin{align} 
\max_{i=0,\ldots,N-1} \E\big| Y_{t_i} - \widehat\Uc_{i}^{(2)}(X_{t_{i}}) \big|^2 & \leq \; C  \E \big|g(\Xc_{T}) - g(X_T) \big|^2 + C |\pi| + C \eps^Z(\pi)  \\
& \hspace{6mm}  + \;  C N \sum_{i=0}^{N-1}  \E \big| \widehat\Vc_{t_i}  -  \widehat\Uc_{i}^{(2)}(X_{t_{i}}) \big|^2.
\end{align}
Plugging \eqref{estimU2} into this last inequality, together with \eqref{convgamma},   gives the required estimation \eqref{eq:theo2_scheme2} for the $Y$-component. 
Finally, by following the same arguments as in  Step 5 in the proof of  \eqref{theo:scheme1_1}, we obtain the estimation \eqref{eq:theo2_scheme2} for the $Z$-component. 
\ep

\vspace{3mm}

\begin{Remark}
\label{remNN_approxError}
{\rm The universal approximation theorem (II) \cite{hor90universal} is valid on compact sets, and one cannot conclude {\it a priori} that the error of network approximation $\eps_i^{\Nc,m}$ converge to zero as $m$ goes to infinity. Instead, we have to proceed into two steps:
\begin{itemize}
    \item[(i)] Localize the error by considering
\begin{align}
\eps_i^{\Nc,m,K} & := \; \inf_{\theta\in\Theta_m^{\gamma}}
\E \big[ \Delta_i(X_{t_i};\theta) 1_{|X_{t_i}| \leq K}   \big], 
\end{align}
where we set $\Delta_i(x;\theta)$ $:=$ $|\hat v_i(x)-\Uc_i(x;\theta)|^2$ + $\Delta t_i \big|\sigma\trans(t_i,x) \big(  D_x \hat v_i(x) - D_x \Uc_i(x;\theta) \big) \big|^2$.
\item[(ii)] Consider an increasing family of neural networks $\Theta_m^{\gamma^{N-1}}$ $\subset$ $\ldots$ $\subset$ $\Theta_m^{\gamma^{i}}$ $\subset$ $\ldots$ 
$\subset$ $\Theta_m^{\gamma^0}$ on which to minimize the approximation errors by backward induction at times $t_i$, $i$ $=$ $N-1,\ldots,0$, and where, $\gamma_m^i$ is defined by  
$$ \gamma^i_m:=\gamma_{\varphi^{N-1-i}(m)},
$$ 
with  $\varphi: \mathbb{N} \to \mathbb{N}$ an increasing function, and where we use the notation $\varphi^k:=\varphi \circ ... \circ \varphi$ (composition $k$ times).

The localized approximation error at time $t_i$, for $0\leq i \leq N-1$, should then be rewritten as
\begin{align}
\eps_{i,N}^{\Nc,m,K} & := \; \inf_{\theta\in\Theta_m^{\gamma^i}}
\E \big[ \Delta_i(X_{t_i};\theta) 1_{|X_{t_i}| \leq K}   \big], \nonumber
\end{align}
and the non-localized one as
\begin{align}
\eps_{i,N}^{\Nc,m} & := \; \inf_{\theta\in\Theta_m^{\gamma^i}}
\E \big[ \Delta_i(X_{t_i};\theta)\big]. \nonumber
\end{align}
\end{itemize}
Note that $\eps_{i,N}^{\Nc,m,K}$ converges to zero, as $m$ goes to infinity, for any $K$ $>$ $0$, as claimed by the universal approximation theorem (II) \cite{hor90universal}. On the other hand, from the expressions of $\hat v_i$, $D_x\hat v_i$ in the above proof of Theorem \ref{theoconv2}, we see under {\bf (H1)}-{\bf (H2)}, and from \eqref{bounderivU} that for all $x$ $\in$ $\R^d$, $\theta$ $\in$ $\Theta_m^{\gamma^i}$, $i$ $=$ $0,\ldots,N-1$: 
\begin{align}
|\Delta_i(x;\theta)| & \leq \; C(1 +|x|^2)\gamma_{\varphi^{N-1}(m)}^4,      
\end{align}
for some positive constant $C$ independent of $m,\pi$. We deduce by Cauchy-Schwarz and Chebyshev's inequalities that for all $K$ $>$ $0$, and 
$\theta$ $\in$ $\Theta_m^{\gamma^i}$, $i$ $=$ $0,\ldots,N-1$, 
\begin{align} \label{deltainter}
\E \big[ \Delta_i(X_{t_i};\theta) 1_{|X_{t_i}| > K} \big] & \leq   
\; \Big\|\Delta_i(X_{t_i};\theta) \Big\|_{_2} \frac{\big\| X_{t_i} \big\|_{_2}}{K}  \; \leq C(1 +|x_0|^3)\frac{\gamma_{\varphi^{N-1}(m)}^4}{K},
\end{align}
where we used \eqref{integX} in the last inequality.  This shows that
\begin{align}
\eps_{i,N}^{\Nc,m} & \leq \;  \eps_{i,N}^{\Nc,m,K} + C \frac{\gamma_{\varphi^{N-1}(m)}^4}{K}, \;\;\; \forall K > 0,    
\end{align}
and thus, in theory, the error $\eps_{i,N}^{\Nc,m}$ can be made arbitrary small by suitable choices of large $m$ and $K$. 
\ep
}
\end{Remark}

\subsection{Convergence of RDBDP}

 In this paragraph, we study the convergence of machine learning schemes for the variational ine\-quality \eqref{eq:IQV}.  
% in the case where the generator $f$ does not depend on $z$. The general case is more difficult to handle and is left for future research. 

We first consider the case  when $f$ does not depend on $z$, so that the component $Y_t$ $=$ $u(t,\Xc_t)$ solution to the reflected BSDE \eqref{RBSDE}  admits 
a Snell envelope representation, and we shall focus on the error on $Y$ by  proposing  an alternative to scheme \eqref{eq:schemeVI}, refereed to as RDBDPbis scheme, which only uses neural network for learning the function $u$: 
 \begin{itemize}
\item  Initialize $\widehat\Uc_N$ $=$ $g$
\item  For $i$ $=$ $N-1,\ldots,0$, given $\widehat\Uc_{i+1}$, use a deep neural network $\Uc_i(.;\theta)$ $\in$ $\Nc\Nc_{d,1,L,m}^\varrho(\R^{N_m})$, 
and compute (by SGD) the minimizer of the expected quadratic loss function
\begin{equation} \label{eq:schemeVIbis}
\left\{ 
\begin{aligned} 
\bar L_i(\theta) & := \E \big| \widehat\Uc_{i+1}(X_{t_{i+1}}) - \Uc_i(X_{t_i};\theta) +  f(t_i,X_{t_i},\Uc_i(X_{t_i};\theta)) \Delta t_i   \big|^2  \\
\theta_i^* & \in {\rm arg}\min_{\theta\in\R^{N_m}} \bar L_i(\theta). 
\end{aligned}
\right.
\end{equation}
Then, update: $\widehat\Uc_i$ $=$ $\max\big[\Uc_i(.;\theta_i^*),g]$. 
\end{itemize}

Let us also define  from the scheme \eqref{eq:schemeVIbis} 
\begin{equation} \label{defV3}
\left\{
 \begin{aligned} 
\tilde\Vc_{t_i} & :=   \; \E_i \big[ \widehat\Uc_{i+1}^{}(X_{t_{i+1}}) \big]  + f(t_i,X_{t_i},\tilde\Vc_{t_i}) \Delta t_i \; = \; \tilde v_i(X_{t_i}), \\
\widehat\Vc_{t_i} & := \;  \max[\tilde\Vc_{t_i} ;  g(X_{t_i}) ],  \;\;\;  i =0,\ldots,N-1.   
\end{aligned}
\right.
\end{equation}

Our next  result  gives an error estimate of the  scheme   \eqref{eq:schemeVIbis} 
in terms of the $L^2$-approximation errors of  $\tilde v_i$  by neural networks $\Uc_i$, $i=0,\ldots,N-1$,  and defined as
\begin{align}
\tilde\eps_i^{\Nc} &  := \; \inf_{\theta\in\R^{N_m}}  \E \big|\tilde v_i(X_{t_i}) - \Uc_i(X_{t_i};\theta) \big|^2. 
\end{align}

\begin{Theorem}   \emph{(Case $f$  independent of $z$:  Consistency of RDBDPbis)}
\label{theo:scheme3bis} Let Assumption {\bf (H1)} hold, with  $g$  Lipschitz. Then, there exists a constant $C>0$, independent of $\pi$, such that
\begin{align}
 \max_{i=0,\ldots,N-1} \big\| Y_{t_i} -  \hat\Uc_i(X_{t_i}) \big\|_{_2}   & \leq \; 
C  \Big(    |\pi|^{\frac{1}{2}}  +  \sum_{i=0}^{N-1} \sqrt{\tilde \eps_i^{\Nc}} \Big),
\label{eq:theo3bis}
\end{align}
where $\|.\|_{_2}$ is the $L^2$-norm on $(\Omega,\Fc,\P)$. 
\end{Theorem}

\begin{Remark} \label{rem:estimRDBDP}
{\rm The estimation \eqref{eq:theo3bis} implies that 
\begin{align}
 \max_{i=0,\ldots,N-1} \E\big| Y_{t_i} -  \hat\Uc_i(X_{t_i}) \big|^2  & \leq \; 
C  \Big(    |\pi|  +  N \sum_{i=0}^{N-1} \tilde \eps_i^{\Nc} \Big),
\label{eq:theo3_scheme3}
\end{align}
which is of the same order than the error estimate in Theorem \ref{theo:scheme1_1} when $g$ is Lipschitz. 
}
\ep 
\end{Remark}

\noindent {\bf Proof.}
Let us  introduce the  discrete-time approximation of the reflected BSDE
\begin{equation} \label{defYpi}
\left\{
\begin{aligned}
Y_{t_N}^\pi & = \; g(X_{t_N})  \\
\tilde Y_{t_i}^\pi &  = \;  \E_i[ Y_{t_{i+1}}^\pi ] +  f(t_i,X_{t_i},\tilde Y_{t_i}^\pi) \Delta t_i \\
Y_{t_i}^\pi &  = \; \max \big[  \tilde Y_{t_i}^\pi  ; g(X_{t_i}) \big] , \;\;\;  i =0,\ldots,N-1. 
\end{aligned}
\right. 
\end{equation}
It is known, see \cite{balpag03}, \cite{bouchard2004discrete} that 
\begin{align} \label{estimRBSDE}
 \max_{i=0,\ldots,N-1} \big\| Y_{t_i} -  Y_{t_i}^\pi  \big\|_{_2}   & = \; O(|\pi|^{\frac{1}{2}}). 
\end{align}
 %By using the fact that $|\max(a,c) - \max(b,c)|$ $\leq$ $|a-b|$, we have for all $i$ $=$ $0,\ldots,N-1$, 
Fix $i$ $=$ $0,\ldots,N-1$.  From \eqref{defV3}, \eqref{defYpi}, we have 
 \begin{align}
 | \tilde Y_{t_i}^\pi - \tilde\Vc_{t_i} | & \leq \; \E_i \big|  Y_{t_{i+1}}^\pi -  \widehat\Uc_{i+1}^{}(X_{t_{i+1}})  \big| 
 + \Delta t_i \big|  f(t_i,X_{t_i},\tilde Y_{t_i}^\pi) - f(t_i,X_{t_i},\tilde\Vc_{t_i}) \big| \\
 & \leq \;  \E_i \big|  Y_{t_{i+1}}^\pi -  \widehat\Uc_{i+1}^{}(X_{t_{i+1}})  \big|  + [f]_{_L} \Delta t_i  | \tilde Y_{t_i}^\pi - \tilde\Vc_{t_i} |,
 \end{align}
from the Lipschitz condition on $f$ in {\bf (H1)}, and then for $|\pi|$ small enough
\begin{align}
\big\|  \tilde Y_{t_i}^\pi - \tilde\Vc_{t_i} \big\|_{_2} & \leq \; (1 + C |\pi|) \big\|  Y_{t_{i+1}}^\pi -  \widehat\Uc_{i+1}^{}(X_{t_{i+1}})  \big\|_{_2}. 
\end{align}
By Minkowski inequality, this  yields for all $\theta$ 
\begin{align} \label{YVxi}
\big\|  \tilde Y_{t_i}^\pi -  \Uc_{i}^{}(X_{t_{i}};\theta)  \big\|_{_2} & \leq \; (1 + C |\pi|) \big\|  Y_{t_{i+1}}^\pi -  \widehat\Uc_{i+1}^{}(X_{t_{i+1}})  \big\|_{_2} + 
\big\| \tilde\Vc_{t_i}  -    \Uc_{i}^{}(X_{t_{i}};\theta)   \big\|_{_2}. 
\end{align}
 
On the other hand, by the martingale representation theorem, there exists an $\R^d$-valued square integrable process $(\tilde Z_t)_t$ such that 
\begin{align} \label{FBSDE3}
\widehat\Uc_{i+1}^{}(X_{t_{i+1}}) & = \;  \tilde\Vc_{t_i} -  
f(t_i,X_{t_i},\tilde\Vc_{t_i}) \Delta t_i  + \int_{t_i}^{t_{i+1}}  \tilde Z_s\trans \diff W_s,  
\end{align}
and the expected squared loss function of the DBDP3 scheme can be  written as 
\begin{align} \label{LtildeL3}
\bar L_i^{}(\theta) &= \;  \tilde L_i(\theta) + \E \Big[ \int_{t_i}^{t_{i+1}} \big| \tilde Z_t  \big|^2 \diff t \Big]
\end{align}
with 
\begin{align}
\sqrt{ \tilde L_i(\theta) } & := \;  \Big\| \tilde\Vc_{t_i}  - \Uc_i(X_{t_i};\theta) + 
\big( f(t_i,X_{t_i},\Uc_i(X_{t_i};\theta)) - f(t_i,X_{t_i},\tilde\Vc_{t_i}) \big) \Delta t_i \Big\|_{_2}.  
%\\ & \hspace{2cm}  + \;  \Delta t_i \E \big| \overline{{\tilde Z_{t_i}}} - \Zc_i(X_{t_i};\eta) \big|^2,
\end{align}
From the Lipschitz condition on $f$, and  by Minkowski inequality,  we have for all $\theta$ 
\begin{eqnarray*}
(1 - [f]_{_L} \Delta t_i)  \big\| \tilde\Vc_{t_i}  - \Uc_i(X_{t_i};\theta) \big\|_{_2}   & \leq & 
\sqrt{\tilde L_i(\theta)}  \; \leq \;  (1 + [f]_{_L} \Delta t_i) \big\| \tilde \Vc_{t_i}  - \Uc_i(X_{t_i};\theta) \big\|_{_2}.   
\end{eqnarray*}
Take now $\theta_i^*$  $\in$ ${\rm arg}\min_\theta \bar L_i(\theta)$ $=$ ${\rm arg}\min_\theta \tilde L_i(\theta)$.  Then, from the above relations, we have
\begin{eqnarray*}
(1 - [f]_{_L} \Delta t_i) \big\| \tilde\Vc_{t_i}  - \Uc_i(X_{t_i};\theta_i^*) \big\|_{_2}   & \leq &  
(1 + [f]_{_L} \Delta t_i)   \big\| \tilde \Vc_{t_i}  - \Uc_i(X_{t_i};\theta) \big\|_{_2}, 
\end{eqnarray*}
for all $\theta$, and so 
\begin{align} \label{VUxi}
\big\| \tilde\Vc_{t_i}  - \Uc_i(X_{t_i};\xi_i^*) \big\|_{_2} & \leq \; (1 + C|\pi|) \sqrt{  \tilde\eps_i^{\Nc} }.   
\end{align}
%Recalling that  $\widehat\Vc_{t_i}$ $=$ $\max[\tilde\Vc_{t_i};g(X_{t_i})]$, and $\hat\Uc_i(X_{t_i})$ $=$ $\max[\Uc_i(X_{t_i};\xi_i^*); g(X_{t_i})]$, and since 
%$|\max(a,c) - \max(b,c)|$ $\leq$ $|a-b|$, we deduce that
%\begin{align}
%\E \big| \widehat\Vc_{t_i}  - \hat\Uc_i(X_{t_i}) \big|^2 & \leq \; C \big( \eps_i^{NN,\tilde v} + \Delta t_i \eps_i^{NN,\tilde z} \big). 
%\end{align}

By taking $\theta$ $=$ $\theta_i^*$ in \eqref{YVxi}, recalling that $\widehat\Uc_i(X_{t_i})$ $=$ $\max[\Uc_i(X_{t_i};\theta_i^*); g(X_{t_i})]$, $Y_{t_i}^\pi$ $=$ $\max[\tilde Y_{t_i}^\pi;g(X_{t_i})]$, 
and since  $|\max(a,c) - \max(b,c)|$ $\leq$ $|a-b|$,  we obtain by using \eqref{VUxi}
\begin{align}
\big\|  Y_{t_i}^\pi -  \widehat \Uc_{i}^{}(X_{t_{i}})  \big\|_{_2} & \leq \; (1 + C |\pi|) \Big( \big\|  Y_{t_{i+1}}^\pi -  \widehat\Uc_{i+1}^{}(X_{t_{i+1}})  \big\|_{_2} + 
\sqrt{ \tilde\eps_i^{\Nc} }  \Big). 
\end{align}
By induction, this yields
\begin{align}
\max_{i=0,\ldots,N-1} \big\|  Y_{t_i}^\pi -  \widehat \Uc_{i}^{}(X_{t_{i}})  \big\|_{_2} & \leq \; C \sum_{i=0}^{N-1} \sqrt{ \tilde\eps_i^{\Nc} },   
\end{align}
and we conclude with \eqref{estimRBSDE}. 
\ep

\vspace{5mm}

We finally turn to the general case when $f$ may depend on $z$, and study  the convergence of the RDBDP scheme \eqref{eq:schemeVI} towards the variational inequality \eqref{eq:IQV} related to the solution $(Y,Z)$ of the reflected BSDE \eqref{RBSDE}  by showing an error estimate for 
\begin{align}
\Ec\big[(\widehat\Uc^{},\widehat\Zc^{}),(Y,Z)\big] & := \; \max_{i=0,\ldots,N-1} \E \big|Y_{t_i}- \widehat\Uc_i^{}(X_{t_i})\big|^2 + \E \bigg[ \sum_{i=0}^{N-1} \int_{t_i}^{t_{i+1}} 
\big| Z_t - \widehat\Zc_i^{}(X_{t_i}) \big|^2 dt \bigg]. 
\end{align}

Let us  define  from the scheme \eqref{eq:schemeVI} 
\begin{equation} \label{defVZ3}
\left\{
 \begin{aligned} 
\tilde\Vc_{t_i} & :=   \; \E_i \big[ \widehat\Uc_{i+1}(X_{t_{i+1}}) \big]  + f(t_i,X_{t_i},\tilde\Vc_{t_i},\overline{{\tilde Z_{t_i}}}) \Delta t_i \; = \; \tilde v_i(X_{t_i}), \\
\overline{{\tilde Z_{t_i}}} & := \;   \frac{1}{\Delta t_i} \E_i\left[ \widehat\Uc_{i+1}(X_{t_{i+1}}) \Delta W_{t_i} \right] \; = \; \tilde z_i(X_{t_i}), \\
\widehat\Vc_{t_i} & := \;  \max[\tilde\Vc_{t_i} ;  g(X_{t_i}) ],  \;\;\;  i =0,\ldots,N-1.   
\end{aligned}
\right.
\end{equation}

Our final  main result  gives an error estimate of the  RDBDP scheme  in terms of the $L^2$-approximation errors of $\tilde v_i$ and $\tilde z_i$ by neural networks $\Uc_i$ and $\Zc_i$, $i=0,\ldots,N-1$, assumed to be independent (see Remark \ref{remNN}), and defined as
\begin{align}
\eps_i^{\Nc,\tilde v} \; := \; \inf_{\xi} \E \big|\tilde v_i(X_{t_i}) - \Uc_i(X_{t_i};\xi) \big|^2, 
\hspace{7mm}  
 \eps_i^{\Nc,\tilde z} \; := \; \inf_{\eta} \E\big|\tilde z_i(X_{t_i}) -  \Zc_i(X_{t_i};\eta) \big|^2. 
\end{align}

The result is obtained under one of  the following additional assumptions 

\vspace{2mm}

\noindent {\bf (H3)} $g$ is $C^1$, and $g$, $D_x g$ are Lipschitz. 

or 

\noindent {\bf (H4)} $\sigma$ is $C^1$, with $\sigma$, $D_x \sigma$ both Lipschitz, and $g$ is $C^2$, with $g$, $D_x g$, $D_x^2 g$  all Lipschitz.

\begin{Theorem} \emph{(Consistency of RDBDP)}
\label{theo:scheme3}
Let Assumption {\bf (H1)} hold. There exists a constant $C>0$, independent of $\pi$, such that
\begin{align}
\Ec\big[(\widehat\Uc^{},\widehat\Zc^{}),(Y,Z)\big] & \leq \; 
C  \Big(  \eps(\pi) +   \sum_{i=0}^{N-1} \big(N \eps_i^{\Nc,\tilde v}  +   \eps_i^{\Nc,\tilde z}\big) \Big), 
\label{eq:theo3}
\end{align}
with $\eps(\pi)$ $=$ $O(|\pi|^{\frac{1}{2}})$ under {\bf (H3)}, and $\eps(\pi)$ $=$ $O(|\pi|)$ under {\bf (H4)}. 
\end{Theorem}
\noindent {\bf Proof.}
Let us  introduce the  discrete-time approximation of the reflected BSDE
\begin{equation} \label{defYZpi}
\left\{
\begin{aligned}
Y_{t_N}^\pi & = \; g(X_{t_N})  \\
Z_{t_i}^\pi &= \; \frac{1}{\Delta t_i} \E_i \big[ Y_{t_{i+1}}^\pi \Delta W_{t_i} \big], \\
\tilde Y_{t_i}^\pi &  = \;  \E_i[ Y_{t_{i+1}}^\pi ] +  f(t_i,X_{t_i},\tilde Y_{t_i}^\pi,Z_{t_i}^\pi) \Delta t_i \\
Y_{t_i}^\pi &  = \; \max \big[  \tilde Y_{t_i}^\pi  ; g(X_{t_i}) \big] , \;\;\;  i =0,\ldots,N-1. 
\end{aligned}
\right. 
\end{equation}
It is known from \cite{boucha08} that 
\begin{equation} \label{estimBC}
\left\{    
\begin{aligned} 
\max_{i=0,\ldots,N-1} \E \big| Y_{t_i} - Y_{t_i}^\pi  \big|^2  & = \; \eps(\pi) \\
\E \bigg[ \sum_{i=0}^{N-1} \int_{t_i}^{t_{i+1}}  \big| Z_t -  Z_{t_i}^\pi \big|^2 dt \bigg] & = \; O(|\pi|^{\frac{1}{2}}),
\end{aligned}
\right.
\end{equation}
with $\eps(\pi)$ $=$ $O(|\pi|^{\frac{1}{2}})$ under {\bf (H3)}, and $\eps(\pi)$ $=$ $O(|\pi|)$ under {\bf (H4)}. 

\vspace{1mm}

Fix $i$ $=$ $0,\ldots,N-1$. By writing that 
\begin{align} 
\tilde Y_{t_i}^\pi  - \tilde\Vc_{t_i} & =  \E_i\big[ Y_{t_{i+1}} - \widehat\Uc_{i+1}^{}(X_{t_{i+1}}) \big] 
+ \Delta t_i \Big(  f(t_i,X_{t_i},\tilde Y_{t_i}^\pi,Z_{t_i}^\pi) -  f(t_i,X_{t_i},\tilde\Vc_{t_i},\overline{{\tilde Z_{t_i}}}) \Big),
\end{align}
and proceeding similarly as in Step 1 in  the proof of  Theorem \ref{theo:scheme1_1}, we have by Young inequality and Lipschitz condition on $f$
\begin{align}
\E\big| \tilde Y_{t_i}^\pi - \tilde\Vc_{t_i} \big|^2 & \leq \; 
(1 +\gamma\Delta t_i) \E \Big|  \E_i\big[ Y_{t_{i+1}}^\pi - \widehat\Uc_{i+1}^{}(X_{t_{i+1}}) \big] \Big|^2 \\
& \;\;\; + 2 \frac{[f]^2_{_L}}{\gamma}  \big(1+ \gamma \Delta t_i\big) 
\Big\{  \Delta t_i \E \big| \tilde Y_{t_i}^\pi  - \tilde\Vc_{t_i} \big|^2   +  \Delta t_i \E \big| Z_{t_i}^\pi - \overline{{\tilde Z_{t_i}}} \big|^2   \Big\}.   \label{Ypiinter}
\end{align} 
From \eqref{defVZ3}, \eqref{defYZpi}, Cauchy-Schwarz inequality, and law of iterated conditional expectations, we have similarly as  in Step 1 in  the proof of  Theorem \ref{theo:scheme1_1}:
\begin{align}
 \Delta t_i \E \big| Z_{t_i}^\pi - \overline{{\tilde Z_{t_i}}} \big|^2& \leq \;  2 d \Big( \E \big|Y_{t_{i+1}}^\pi- \widehat\Uc_{i+1}^{}(X_{t_{i+1}}) \big|^2 - 
\E \Big| \E_i\big[ Y_{t_{i+1}}^\pi  - \widehat\Uc_{i+1}^{}(X_{t_{i+1}})\big]   \Big|^2  \Big).  
\end{align}
Then,  by  plugging into \eqref{Ypiinter} and choosing $\gamma$ $=$ $4 d [f]^2_{_L}$, we have   for $|\pi|$ small enough: 
\begin{align}
\E\big| \tilde Y_{t_i}^\pi - \tilde\Vc_{t_i} \big|^2 & \leq \; (1 + C |\pi|) \E \big|Y_{t_{i+1}}^\pi - \widehat\Uc_{i+1}^{}(X_{t_{i+1}}) \big|^2.   
\end{align}
 
Next, by using Young inequality as in Step 2 in  the proof of  Theorem \ref{theo:scheme1_1}, we obtain for all $\theta$ $=$ $(\xi,\zeta)$:
\begin{align} \label{tildeYU}
\hspace{-5mm} \E\big| \tilde Y_{t_i}^\pi - \Uc_{i}(X_{t_i};\xi) \big|^2 & \leq \; (1 + C |\pi|) \E \big|Y_{t_{i+1}}^\pi - \widehat\Uc_{i+1}^{}(X_{t_{i+1}}) \big|^2 
+  CN   \E\big| \tilde\Vc_{t_i} - \Uc_{i}(X_{t_i};\xi) \big|^2. 
\end{align}

On the other hand, by the martingale representation theorem, there exists an $\R^d$-valued square integrable process $(\tilde Z_t)_t$ such that 
\begin{align} \label{RFBSDE}
\widehat\Uc_{i+1}^{}(X_{t_{i+1}}) & = \;  \tilde\Vc_{t_i} -  
f(t_i,X_{t_i},\tilde\Vc_{t_i},\overline{{\tilde Z_{t_i}}}) \Delta t_i  + \int_{t_i}^{t_{i+1}}  \tilde Z_s\trans \diff W_s,  
\end{align}
and the expected squared loss function of the RDBDP scheme can be  written as 
\begin{align} 
\hat L_i(\theta) &= \;  \tilde L_i(\theta) + \E \Big[ \int_{t_i}^{t_{i+1}} \big| \tilde Z_t - \overline{{\tilde Z_{t_i}}} \big|^2 \diff t \Big], 
\end{align}
where we notice  by It\^o isometry that $\overline{{\tilde Z_{t_i}}}$ $=$ $\frac{1}{\Delta t_i}\E_i\Big[ \int_{t_i}^{t_{i+1}}  \tilde Z_t dt \Big]$, and 
\begin{align}
 \tilde L_i(\theta)  & := \; \E \Big| \tilde\Vc_{t_i}  - \Uc_i(X_{t_i};\xi) + 
\big( f(t_i,X_{t_i},\Uc_i(X_{t_i};\xi), \Zc_i(X_{t_i};\eta)) -f(t_i,X_{t_i},\tilde\Vc_{t_i},\overline{{\tilde Z_{t_i}}}) \big) \Delta t_i \Big|^2 \\ 
& \hspace{.6cm}  + \;  \Delta t_i \E \big| \overline{{\tilde Z_{t_i}}} - \Zc_i(X_{t_i};\eta) \big|^2. 
\end{align}
By the same arguments as in Step 3 in  the proof of  Theorem \ref{theo:scheme1_1},  using Lipschitz condition on $f$ and Young inequality, we show that for all 
$\theta$ $=$ $(\xi,\eta)$ 
\begin{align}
(1 - C \Delta t_i) \E \big| \tilde\Vc_{t_i}  - \Uc_i(X_{t_i};\xi) \big|^2 +   \frac{\Delta t_i}{2} \E \big| \overline{{\tilde Z_{t_i}}} - \Zc_i(X_{t_i};\eta) \big|^2 \\
& \hspace{-7cm}\leq \tilde L_i(\theta)  \; \leq  \; (1 + C \Delta t_i) \E \big| \tilde\Vc_{t_i}  - \Uc_i(X_{t_i};\xi) \big|^2 + C \Delta t_i \E \big| \overline{{\tilde Z_{t_i}}} - \Zc_i(X_{t_i};\eta) \big|^2. 
\end{align}
By taking $\theta_i^*$ $=$ $(\xi_i^*,\eta_i^*)$ $\in$ ${\rm arg}\min_\theta \hat L_i^{}(\theta)$ $=$ ${\rm arg}\min_\theta \tilde L_i^{}(\theta)$, it follows that for $|\pi|$ small enough 
\begin{align} 
\E \big| \tilde\Vc_{t_i}  - \Uc_i^{}(X_{t_i};\xi_i^*) \big|^2 +   \Delta t_i \E \big| \overline{{\tilde Z_{t_i}}} - \Zc_i^{}(X_{t_i};\eta_i^*) \big|^2
& \leq \; C \eps_i^{\Nc,\tilde v} + C \Delta t_i \eps_i^{\Nc,\tilde z}. 
\end{align}
By plugging into \eqref{tildeYU}, recalling that $\widehat\Uc_i(X_{t_i})$ $=$ $\max[\Uc_i(X_{t_i};\xi_i^*); g(X_{t_i})]$, $Y_{t_i}^\pi$ $=$ $\max[\tilde Y_{t_i}^\pi;g(X_{t_i})]$, 
and since  $|\max(a,c) - \max(b,c)|$ $\leq$ $|a-b|$, we obtain   
\begin{align} 
\E\big|  Y_{t_i}^\pi - \widehat\Uc_{i}(X_{t_i}) \big|^2 & \leq \; (1 + C |\pi|) \E \big|Y_{t_{i+1}}^\pi - \widehat\Uc_{i+1}^{}(X_{t_{i+1}}) \big|^2 
+  CN   \big(  \eps_i^{\Nc,\tilde v} +  \Delta t_i \eps_i^{\Nc,\tilde z} \big), 
\end{align} 
and then by induction
\begin{align}
\max_{i=0,\ldots,N-1} \E\big|  Y_{t_i}^\pi - \widehat\Uc_{i}(X_{t_i}) \big|^2 & \leq \; C \sum_{i=0}^{N-1} \big(  N \eps_i^{\Nc,\tilde v} +   \eps_i^{\Nc,\tilde z} \big). 
\end{align} 
Combining with \eqref{estimBC}, this proves the error estimate \eqref{eq:theo3} for the $Y$-component. The  error estimate \eqref{eq:theo3} for the $Z$-component is proved along the same arguments as in Step 5 in  the proof of  Theorem \ref{theo:scheme1_1}, and is omitted here.  
\ep

\section{Numerical results}
\label{secnum}

In the first two subsections, we compare our schemes DBDP1 \eqref{eq:scheme1}, DBDP2 \eqref{eq:scheme2} and the scheme proposed by \cite{han2017overcoming} on some examples of PDEs and BSDEs.

We first test our algorithms on some PDEs with bounded solutions and quite a simple structure (see section \ref{sec:numeric1}), and then try to solve some PDEs with unbounded solutions and more complex structures (see section \ref{sec:numeric2}).
Our goal is to emphasize that solutions with simple structure easily represented by a neural network can be evaluated by our method even in very high-dimension, whereas the solution with complex structure can only be evaluated in moderate dimension.

Finally, we apply the scheme described in section \ref{sec:varIneq} to an American option problem and show its accuracy in high dimension (see section \ref{sec:numeric3}).

If not specified, we use in the sequel a fully connected feedforward network with two hidden layers, and $d+10$ neurons on each hidden layer, to implement our schemes \eqref{eq:scheme1} and \eqref{eq:scheme2}. We  choose tanh as activation function for the hidden layers in order to avoid some explosion while calculating the numerical gradient $Z$ in scheme \eqref{eq:scheme2} and choose identity function as activation function for the output layer.  We renormalize the data before entering the network. We use Adam Optimizer, implemented in TensorFlow and mini-batch with $1000$ trajectories for the stochastic gradient descent. 

\subsection{PDEs with bounded solution and simple structure}
\label{sec:numeric1}

We begin with a simple example in dimension one. It is not hard to find test cases where the scheme proposed in \cite{han2017overcoming} fails even in dimension one. In fact the latter scheme works well for small maturities and  with a starting point close to the solution.

It is always interesting to start by testing schemes in dimension one as one can ea\-sily compare graphically the numerical results to the theoretical solution. Then we take some examples in higher dimensions and show that our method seems to work well when the dimension increases higher.

\subsubsection{An example in 1D}

We take the following parameters for the BSDE problem defined by \eqref{eq:SDE} and \eqref{eqBSDE}:
\begin{equation}
	\sigma = 1 , \; \mu = 0.2 ,  \; T = 2 , \; d=1,
	\label{coeff:pb_simple}
\end{equation} 
\begin{equation}
	\begin{array}{rclrclrclrcl}
	f(t,x,y,z)  &= &  (\cos(x) (e^{\frac{T-t}{2}}+ \frac{\sigma^2}{2}) + \mu \sin(x)) e^{\frac{T-t}{2}} - \frac{1}{2} \left( \sin(x)\cos(x) e^{T-t} \right)^2 + 
	\frac{1}{2}(yz)^2& & & & & & & & & \\
	g(x)&= &\cos(x). & & & & & &  & 
	\end{array}
\end{equation}
for which, the explicit analytic solution is equal to $u(t,x) =e^{\frac{T-t}{2}} \cos(x)$.

We want to estimate the solution $u$ and its gradient $D_x u$ from our schemes. This example is interesting, because  with $T=1$, the method proposed in \cite{han2017overcoming}, initializing  $u(0,.)$ as the solution of the associated linear problem associated ($f=0$) and randomly initializing $D_x u(0,.)$ works very well. However, for $T=2$, the method in \cite{han2017overcoming} always fails on our test whatever the choice of the  initialization: the algorithm is either trapped in a local minimum when the initial learning rate associated to the gradient method is too small or explodes when the learning rate is taken higher.
This numerical failure is not dependent on the considered network: using some LSTM networks as in \cite{chan2018machine} gives the same result.

\begin{figure}[h!]
\begin{minipage}[b]{0.49\linewidth}
  \centering
 \includegraphics[width=\textwidth]{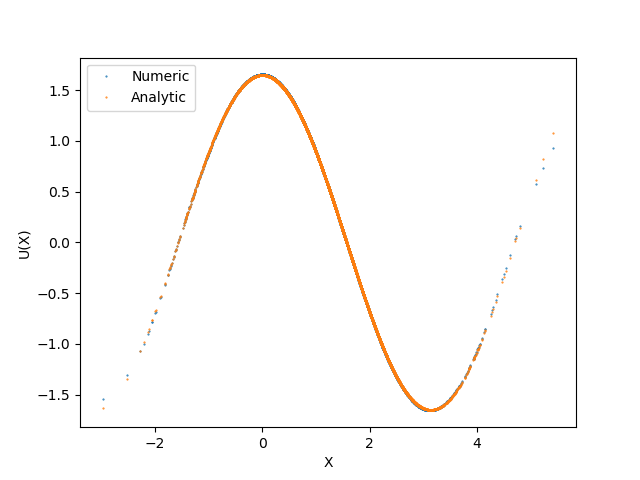}
 \caption*{$u(t,.)$ and its estimate at time $t=1.$}
 \end{minipage}
\begin{minipage}[b]{0.49\linewidth}
  \centering
 \includegraphics[width=\textwidth]{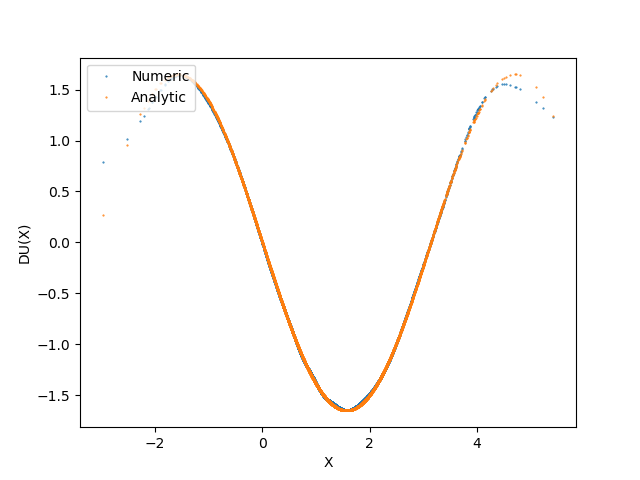}
 \caption*{$Z$ and its estimate at time $t=1.$ }
 \end{minipage}
 \begin{minipage}[b]{0.49\linewidth}
  \centering
 \includegraphics[width=\textwidth]{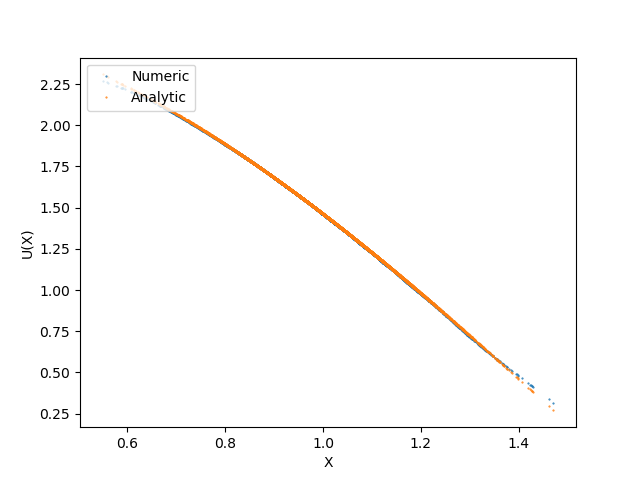}
 \caption*{$u(t,.)$ and its estimate at time $t=0.0091$.}
 \end{minipage}
 \begin{minipage}[b]{0.49\linewidth}
  \centering
 \includegraphics[width=\textwidth]{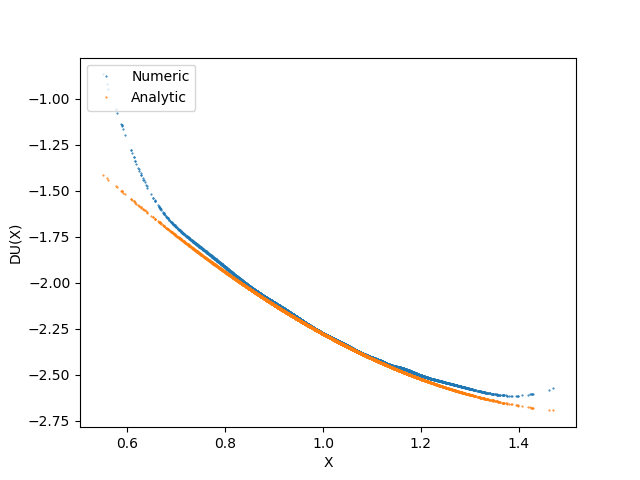}
 \caption*{$Z$ and its estimate at time $t=0.0091$.}
 \end{minipage}
 \caption{\label{fig:quentinCase1DS1} Estimates of $u$ and $Z$ using DBDP1. We took the parameters defined in \eqref{coeff:pb_simple} and set $x_0$ $=$ $1$.  }
 \end{figure}

\begin{figure}[h!]
\begin{minipage}[b]{0.49\linewidth}
  \centering
 \includegraphics[width=\textwidth]{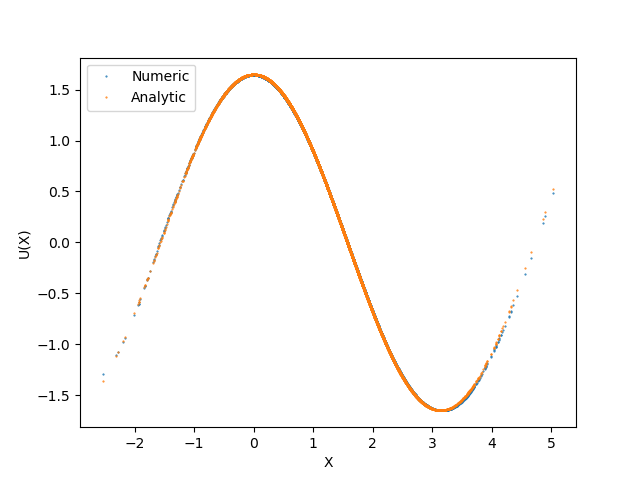}
 \caption*{$u(t,.)$ and its estimate at time $t=1.$}
 \end{minipage}
\begin{minipage}[b]{0.49\linewidth}
  \centering
 \includegraphics[width=\textwidth]{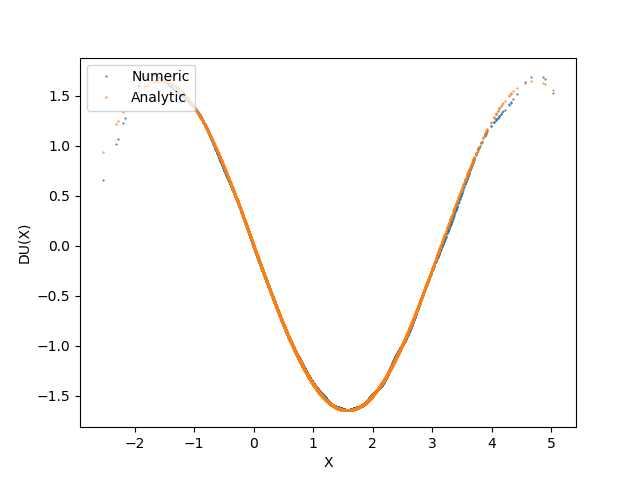}
 \caption*{$Z$ and its estimate at time $t=1.$ }
 \end{minipage}
 \begin{minipage}[b]{0.49\linewidth}
  \centering
 \includegraphics[width=\textwidth]{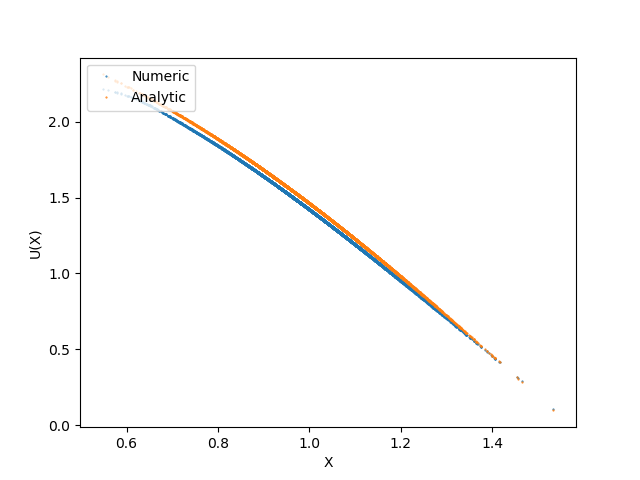}
 \caption*{$u(t,.)$ and its estimate at time $t=0.0091$.}
 \end{minipage}
 \begin{minipage}[b]{0.49\linewidth}
  \centering
 \includegraphics[width=\textwidth]{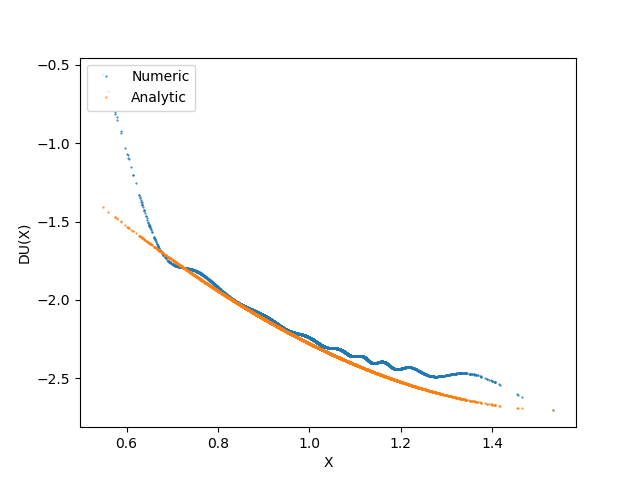}
 \caption*{$Z$ and its estimate at time $t=0.0091$.}
 \end{minipage}
 \caption{\label{fig:quentinCase1DS2} Estimates of $u$ and $Z$ using DBDP2. We took the parameters defined in \eqref{coeff:pb_simple} and set $x_0$ $=$ $1$.}
 \end{figure}

Because of the high non-linearity, we discretize the BSDE using $N$ $=$ $240$ time steps, and implemented hidden layers with $d+10$ $=$ $11$ neurons. Figure \ref{fig:quentinCase1DS1} (resp. Figure \ref{fig:quentinCase1DS2}) depicts the estimated functions $u(t,.)$ and $D_x u(t,.)$ estimated from DBDP1 (resp. DBDP2) scheme.

\begin{table}[H]
\centering

 \begin{tabular}{|c|c|c|}  \hline
 &  Averaged value & Standard deviation \\ \hline
 DBDP1 & 1.46332  & 0.01434  \\ \hline
 DBDP2 & 1.4387982  & 0.01354  \\ \hline
 \end{tabular}
\caption{Estimate of $u(0,x_0)$ where $d$ $=$ $1$ and $x_0$ $=$ $1$. Average and standard deviation observed over 10 independent runs are reported. The theoretical solution is $1.4686938$.}
\label{tab:sol1D}
\end{table}
 \clearpage

\subsubsection{Increasing the dimension}

We extend the example from the previous section to the following $d$-dimensional problem:
\begin{flalign*}
d\geq 1, \quad & \sigma = \frac{1}{\sqrt{d}} \I_d,  \quad  \mu = \frac{0.2}{d} \un_d , \quad T = 1 , 
\end{flalign*}
\begin{equation}
\begin{array}{rcl}
f(t,x,y,z) & = &  (\cos( \bar x) (e^{\frac{T-t}{2}}+ \frac{1}{2}) + 0.2 \sin( \bar x)) e^{\frac{T-t}{2}} - \frac{1}{2} \left( \sin(\bar x)\cos(\bar x) e^{T-t} \right)^2 + \frac{1}{2d}(u(\un_d.z))^2,\\
g(x) & = &\cos(\bar x),
\end{array}
\end{equation}
with $\bar x = \sum_{i=1}^d x_i$.

We take $N$ $=$ $120$ in the Euler scheme, and $d+10$ neurons for each hidden layer. We take $1000$ trajectories in mini batch, use data renormalization, and check the loss convergence every $50$ iterations. For this small maturity,  the  scheme \cite{han2017overcoming} generally converges, and we give the results obtained with the same network and initializing the scheme with the linear solution of the problem.
Results in dimension 5 to 50 are given in Tables \ref{tab:sol5Dsimple}, \ref{tab:sol10Dsimple}, \ref{tab:sol20D} and
\ref{tab:sol50D}. Both schemes \eqref{eq:scheme1} and \eqref{eq:scheme2} work well  with results very close to the solution and close to the results calculated by the scheme \cite{han2017overcoming}.
As the dimension increases, scheme \eqref{eq:scheme1} seems to be the most accurate.

\begin{Remark}
{\rm
In dimension $50$, the  initial learning rate in scheme \cite{han2017overcoming} is taken small in order to avoid a divergence of the method. In fact, running the test 3 times (with 10 runs each time), we observed convergence of the algorithm two times, and in the last test: one of the ten run exploded, and another one clearly converged to a wrong solution.
}
\ep
\end{Remark}

\begin{table}[H]
\centering
\begin{tabular}{|c|c|c|}  \hline
 &  Averaged value & Standard deviation \\ \hline
 DBDP1 & 0.4637038  & 0.004253  \\ \hline
 DBDP2 & 0.46335  & 0.00137 \\ \hline
 Scheme \cite{han2017overcoming}  & 0.46562 &  0.0035\\ \hline
 \end{tabular}
\caption{Estimate of $u(0,x_0)$ where $d$ $=$ $5$ and $x_0$ $=$ $\un_5$. Average and standard deviation observed over 10 independent runs are reported. The theoretical solution is $0.46768$. %Solution in  dimension $d$$=$$5$ at point $x=\un_5$, $t=0$ as the average of 10 runs and standard deviation observed. The solution is $0.46768$.
}
\label{tab:sol5Dsimple} 
\end{table}

%\newpage

\begin{table}[H]
\centering
\begin{tabular}{|c|c|c|}  \hline
 &  Averaged value & Standard deviation \\ \hline
 DBDP1 &  - 1.3895 &  0.00148\\ \hline
 DBDP2 & -1.3913 &  0.000583 \\ \hline
 Scheme \cite{han2017overcoming}  & -1.3880 &  0.00155\\ \hline
 \end{tabular}
\caption{
	Estimate of $u(0,x_0)$ where $d$ $=$ $10$ and $x_0$ $=$ $\un_{10}$. Average and standard deviation observed over 10 independent runs are reported. The theoretical solution is $-1.383395$.
	%Solution in  dimension $10$ at point $x=\un_d $, $t=0$ as the average of  10 runs and standard deviation observed. The solution is $-1.383395$.
}
\label{tab:sol10Dsimple}
\end{table}
\begin{table}[H]
\centering
\begin{tabular}{|c|c|c|}  \hline
 &  Averaged value & Standard deviation \\ \hline
 DBDP1 &  0.6760 & 0.00274 \\ \hline
 DBDP2 &0.67102 & 0.00559 \\ \hline
Scheme \cite{han2017overcoming}  & 0.68686 & 0.002402  \\ \hline
 \end{tabular}
\caption{\label{tab:sol20D} 
		Estimate of $u(0,x_0)$ where $d$ $=$ $20$ and $x_0$ $=$ $\un_{20}$. Average and standard deviation observed over 10 independent runs are reported. The theoretical solution is $0.6728135$.
	%Solution in  dimension $20$ at point $x=\un_d $, $t=0$ as the average of  10 runs and standard deviation observed. The solution is $0.6728135$.
}
\end{table}
\begin{table}[H]
\centering
\begin{tabular}{|c|c|c|}  \hline
 &  Averaged value & Standard deviation \\ \hline
 DBDP1 &  1.5903 & 0.006276 \\ \hline
 DBDP2 & 1.58762& 0.00679 \\ \hline
Scheme \cite{han2017overcoming}  & 1.583023  & 0.0361 \\ \hline
 \end{tabular}
\caption{\label{tab:sol50D} 
Estimate of $u(0,x_0)$ where $d$ $=$ $50$ and $x_0$ $=$ $\un_{50}$. Average and standard deviation observed over 10 independent runs are reported. The theoretical solution is $1.5909$.
	%Solution in  dimension $50$ at point $x=\un_d $, $t=0$ as the average of  10 runs and standard deviation observed. The solution is $1.5909$.
}
\end{table}

\subsection{PDEs with unbounded solution and more complex structure}
\label{sec:numeric2}

\begin{figure}[h!]
\begin{minipage}[b]{0.49\linewidth}
  \centering
 \includegraphics[width=\textwidth]{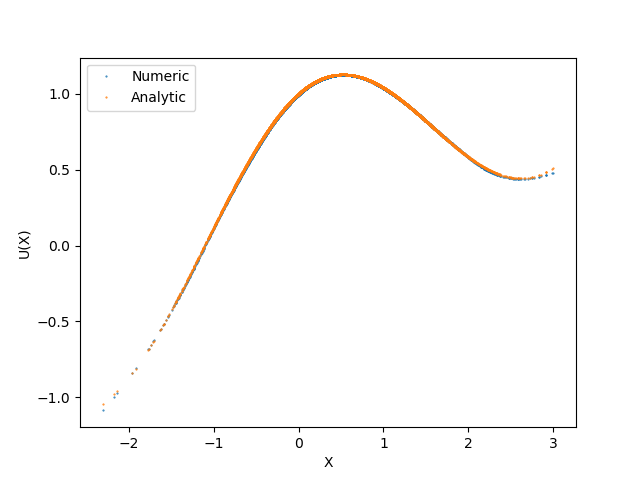}
 \caption*{$u(t,.)$ and its estimate at time $t=0.5$.}
 \end{minipage}
\begin{minipage}[b]{0.49\linewidth}
  \centering
 \includegraphics[width=\textwidth]{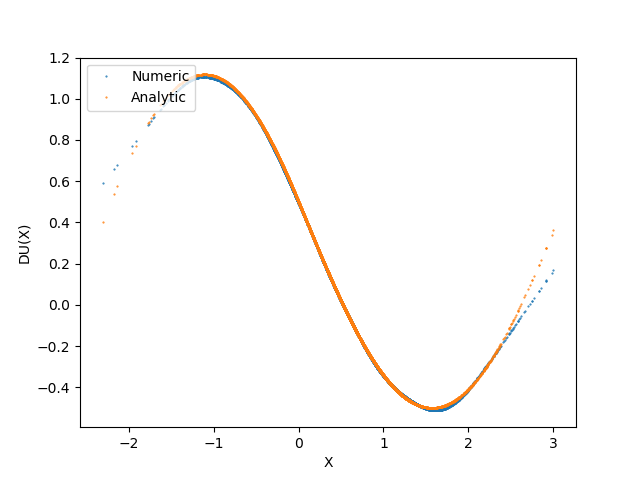}
 \caption*{$Z$ and its estimate at time $t=0.5$ }
 \end{minipage}
 \begin{minipage}[b]{0.49\linewidth}
  \centering
 \includegraphics[width=\textwidth]{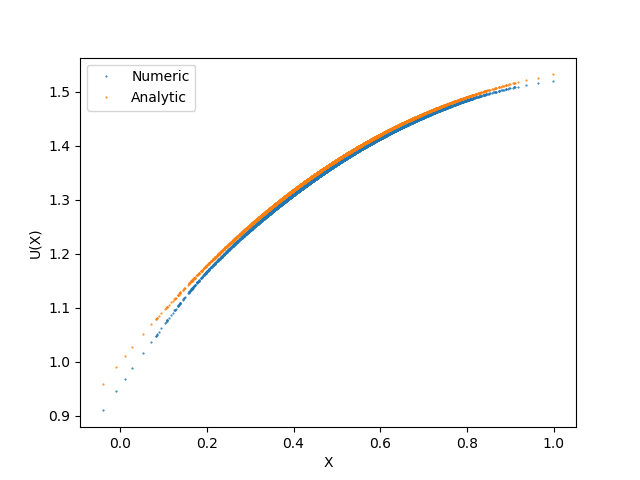}
 \caption*{$u(t,.)$ and its estimate at time $t=0.0085$.}
 \end{minipage}
 \begin{minipage}[b]{0.49\linewidth}
  \centering
 \includegraphics[width=\textwidth]{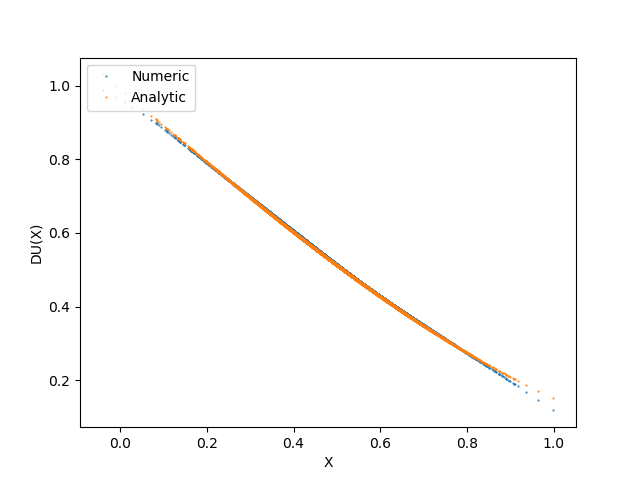}
 \caption*{$Z$ and its estimate at time $t=0.0085$.}
 \end{minipage}
 \caption{\label{fig:huyenCase1DS1} 
 	Estimates of $u$ and $Z$ using DBDP1. We took the parameters defined in \eqref{coeff:pb_complex}, with $d$ $=$ $1$, and set $x_0$ $=$ $0,5$.}
 \end{figure}

\begin{figure}[h!]
\begin{minipage}[b]{0.49\linewidth}
  \centering
 \includegraphics[width=\textwidth]{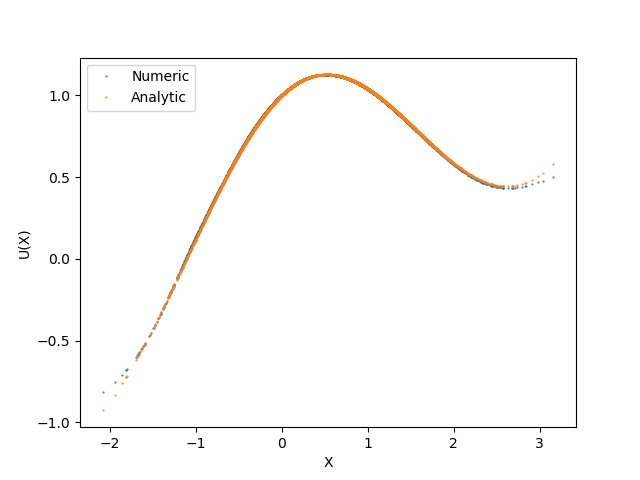}
 \caption*{$u(t,.)$ and its estimate at time $t=0.5$.}
 \end{minipage}
\begin{minipage}[b]{0.49\linewidth}
  \centering
 \includegraphics[width=\textwidth]{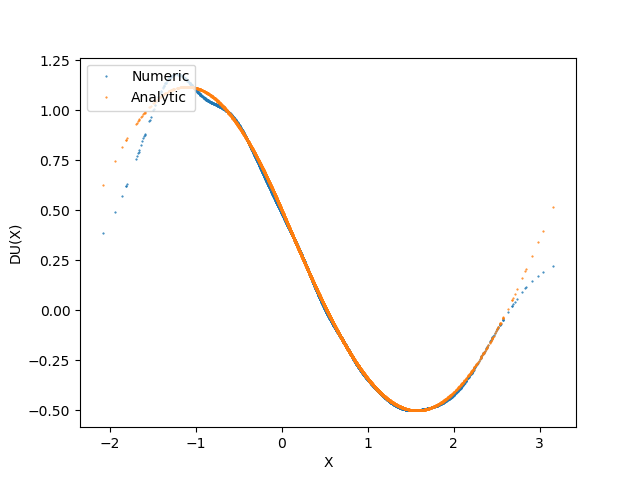}
 \caption*{$Z$ and its estimate at time $t=0.5$ }
 \end{minipage}
 \begin{minipage}[b]{0.49\linewidth}
  \centering
 \includegraphics[width=\textwidth]{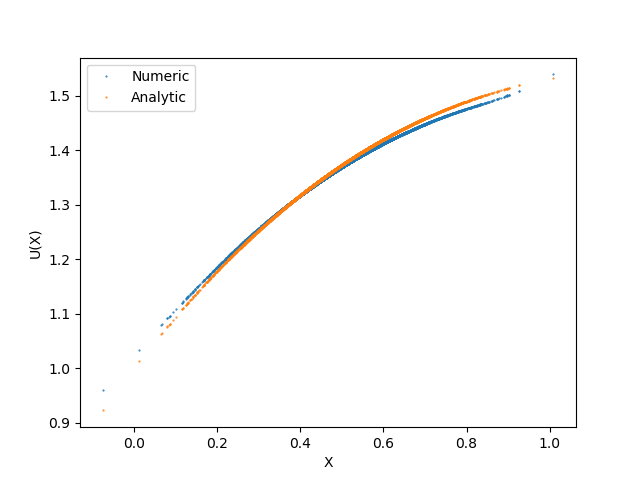}
 \caption*{$u(t,.)$ and its estimate at time $t=0.0085$.}
 \end{minipage}
 \begin{minipage}[b]{0.49\linewidth}
  \centering
 \includegraphics[width=\textwidth]{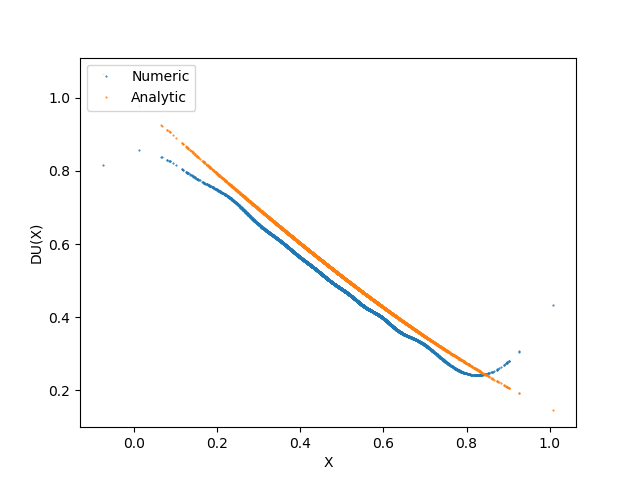}
 \caption*{$Z$ and its estimate at time $t=0.0085$.}
 \end{minipage}
 \caption{\label{fig:huyenCase1DS2}
 	Estimates of $u$ and $Z$ using DBDP2. We took the parameters defined in \eqref{coeff:pb_complex}, with $d$ $=$ $1$, and set $x_0$ $=$ $0.5$.
 	 %Test case 1D with the unbounded solution : results for scheme \eqref{eq:scheme2}, starting point is $0.5$.  
  }
 \end{figure}

In this section with take the following parameters
\begin{flalign}
\sigma  & = \;  \frac{1}{\sqrt{d}} \I_d,  \quad  \mu = 0, \quad T = 1 , \\
f(x,y,z)  & = \;    k(x) + \frac{1}{2 \sqrt{d}} y(\un_d.z) + \frac{ y^2}{2} 
\label{coeff:pb_complex}
\end{flalign}
where the function $k$ is chosen such that the solution to the PDE is equal to
\begin{flalign*}
u(t,x) = \frac{T-t}{d} \sum_{i=1}^d ( \sin(x_i) 1_{x_i<0} + x_i 1_{x_1 \ge 0} ) +  \cos \left(\sum_{i=1}^d i x_i \right).
\end{flalign*}
Notice that the structure of the solution is more complex than in the first example. We aim at evaluating the solution at $x=0.5 \un_d$. We take $120$ time steps for the Euler time discretization and $d+10$ neurons in each hidden layers.
As shown in Figures \ref{fig:huyenCase1DS1} and \ref{fig:huyenCase1DS2} as well as in Table \ref{tab:sol1D_complex}, the three schemes provide accurate and stable results in dimension $d$ $=$ $1$.

 \begin{table}[H]
 \centering
 \begin{tabular}{|c|c|c|}  \hline
 &  Averaged value & Standard deviation \\ \hline
DBDP1 & 1.3720 & 0.00301 \\ \hline
DBDP2 & 1.37357 & 0.0022 \\ \hline
 Scheme \cite{han2017overcoming}  &  1.37238    &  0.00045      \\ \hline
 \end{tabular}
\caption{
	Estimate of $u(0,x_0)$, where $d$ $=$ $1$ and $x_0$ $=$ $0.5$. Average and standard deviation observed over 10 independent runs are reported. The theoretical solution is $1.37758$.
	%Solution in dimension $d$=1 for the  unbounded case at point $x=0.5$, $t=0$ as the average of  10 runs and standard deviation observed. The theoretical solution is $1.37758$.
}
\label{tab:sol1D_complex}
\end{table}

\begin{figure}[h!]
\begin{minipage}[b]{0.49\linewidth}
  \centering
 \includegraphics[width=\textwidth]{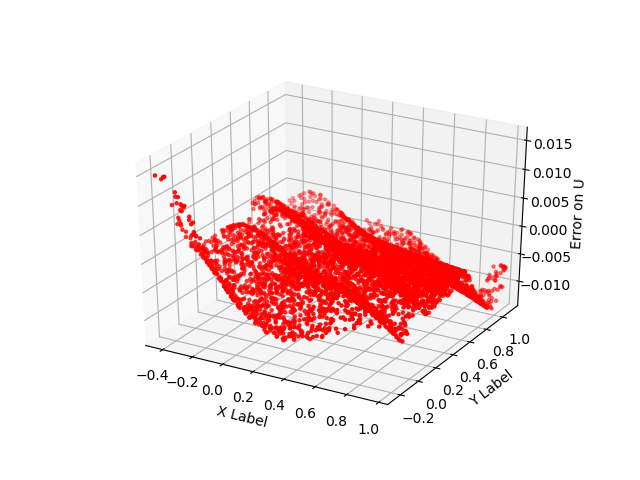}
 \caption*{Error on solution at date $t=0.5$.}
 \end{minipage}
\begin{minipage}[b]{0.49\linewidth}
  \centering
 \includegraphics[width=\textwidth]{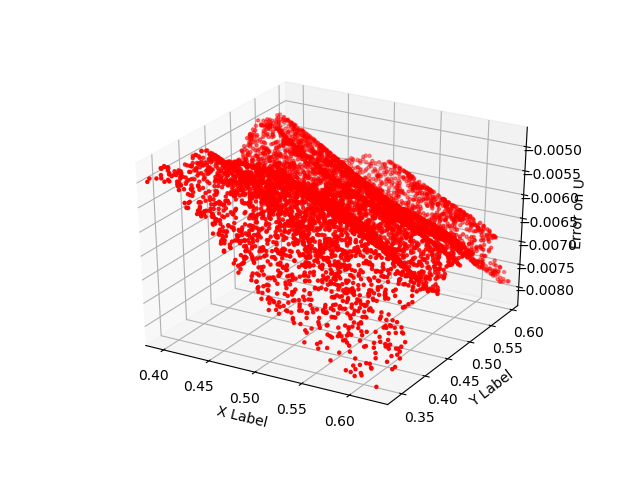}
 \caption*{Error on solution at date $t=0.0085$.}
 \end{minipage}
 \caption{\label{fig:huyenCase2DS1} Algebric error of the estimate of $u$ using DBDP1. We took the parameters in \eqref{coeff:pb_complex} and set $d$ $=$ $2$ and $x_0$ $=$ $0.5  \un_d$.  }
 \end{figure}

\begin{figure}[h!]
\begin{minipage}[b]{0.49\linewidth}
  \centering
 \includegraphics[width=\textwidth]{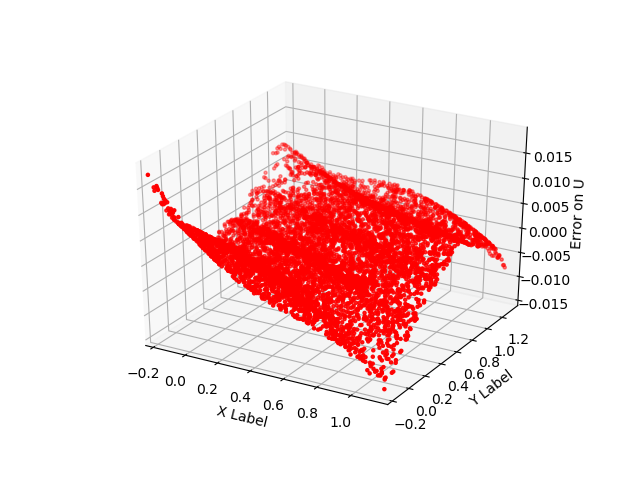}
 \caption*{Error on solution at date $t=0.5$.}
 \end{minipage}
 \begin{minipage}[b]{0.49\linewidth}
  \centering
 \includegraphics[width=\textwidth]{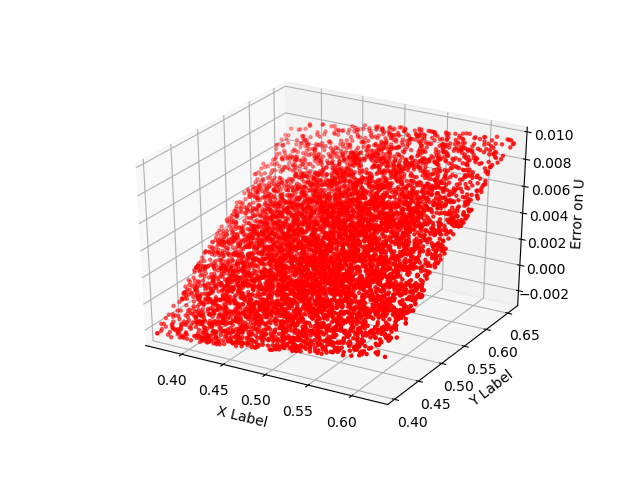}
 \caption*{Error on solution at date $t=0.0085$.}
 \end{minipage}
 \caption{\label{fig:huyenCase2DS2}  Algebric error of the estimate of $u$ using scheme \eqref{eq:scheme2}. We took the parameters in \eqref{coeff:pb_complex} and set $d$ $=$ $2$ and $x_0$ $=$ $0.5  \un_d$.}
 \end{figure}

In dimension 2, the three schemes provide very accurate and stable results, as shown in Figures \ref{fig:huyenCase2DS1} and \ref{fig:huyenCase2DS2}, as well as in Table \ref{tab:sol2D}.

 \begin{table}[H]
 \centering
 \begin{tabular}{|c|c|c|}  \hline
 &  Averaged value & Standard deviation \\ \hline
DBDP1 & 0.5715359 & 0.0038 \\ \hline
DBDP2 & 0.5707974 &  0.00235\\ \hline
 Scheme \cite{han2017overcoming}  &  0.57145   &  0.0006      \\ \hline
 \end{tabular}
\caption{\label{tab:sol2D} 
		Estimate of $u(0,x_0)$, where $d$ $=$ $2$ and $x_0$ $=$ $0.5\un_{2}$. Average and standard deviation observed over 10 independent runs are reported. The theoretical solution is $0.570737$.
		%Solution in  dimension 2 for the unbounded case at point $x=0.5$, $t=0$ as the average of 10 runs and standard deviation observed. The solution is $0.570737$.
	}
\end{table}

Above dimension 3, the scheme \cite{han2017overcoming} always explodes no matter the chosen initial learning rate and the activation function for the hidden layers (among the $\tanh$, ELU, ReLu and sigmoid ones). Besides, taking $3$ or $4$ hidden layers does not improve the results.

We reported the results obtained in dimension $d$ $=$ $5$ and $8$ in Table \ref{tab:sol5D} and \ref{tab:sol8D}. Scheme \eqref{eq:scheme1} seems to work better than scheme \eqref{eq:scheme2} as the dimension increases. Note that the standard deviation increases with the dimension of the problem.

\begin{table}[H]
 \centering
 \begin{tabular}{|c|c|c|}  \hline
 &  Averaged value & Standard deviation \\ \hline
DBDP1 & 0.8666 & 0.013 \\ \hline
DBDP2 & 0.83646 & 0.00453 \\ \hline
 Scheme \cite{han2017overcoming}  &  NC   &     NC   \\ \hline
 \end{tabular}
\caption{\label{tab:sol5D} 
	Estimate of $u(0,x_0)$, where $d$ $=$ $5$ and $x_0$ $=$ $0.5\un_{5}$. Average and standard deviation observed over 10 independent runs are reported. The theoretical solution is $0.87715$.
	%Solution in  dimension 5 for the unbounded case at point $x=0.5$, $t=0$ as the average of 20 runs and standard deviation observed. The solution is $0.87715$. 
}
\end{table}

\begin{table}[H]
 \centering
 \begin{tabular}{|c|c|c|c|}  \hline
 &    Averaged value & Standard deviation \\ \hline
DBDP1 & 1.169441  & 0.02537 \\ \hline
DBDP2 &  1.0758344 & 0.00780 \\ \hline
 Scheme \cite{han2017overcoming}  &  NC   &     NC   \\ \hline
 \end{tabular}
\caption{\label{tab:sol8D} 
		Estimate of $u(0,x_0)$, where $d$ $=$ $8$ and $x_0$ $=$ $0.5\un_{8}$. Average and standard deviation observed over 10 independent runs are reported. The theoretical solution is $1.1603167$.
		%Solution in  dimension 8 for the unbounded case at point $x=0.5$, $t=0$ as the average of  50 runs and standard deviation observed.The solution is $1.1603167$.
	}
\end{table}

When $d\geq10$, schemes \eqref{eq:scheme1} and \eqref{eq:scheme2} both fail at providing correct estimates of the solution, as shown in Table \ref{tab:sol10D}. Increasing the number of layers or neurons does not improve the result.

\begin{table}[H]
 \centering
 \begin{tabular}{|c|c|c|}  \hline
 &  Averaged value & Standard deviation \\ \hline
DBDP1 & -0.3105 &  0.02296 \\ \hline
DBDP2 &  -0.3961&  0.0139\\ \hline
  Scheme \cite{han2017overcoming}  &  NC   &     NC   \\ \hline
 \end{tabular}
\caption{\label{tab:sol10D} 
	Estimate of $u(0,x_0)$, where $d$ $=$ $10$ and $x_0$ $=$ $0.5\un_{10}$. Average and standard deviation observed over 10 independent runs are reported. The theoretical solution is $-0.2148861$.
	%Solution in  dimension 10 for unbounded case at point $x=0.5$, $t=0$ as the average of 50 runs and standard deviation observed. Three hidden layers used. Solution is $-0.2148861$.
}
\end{table}

\subsection{Application to American options}
\label{sec:numeric3}

Consider the stock price $X_t$ $=$ $(X^1_t, \dots , X^d_t)$ of $d$ assets with  the following dynamics under the risk neutral  probability measure:
\begin{equation*}
     dX_t^i  =  r X_t^i dt + \sigma_i X_t^i dW_t^i,
\end{equation*}
where $W_{.} = (W_{.}^1, \dots, W^d_{.})$ is a $d$-dimensional Brownian Motion, $\sigma = ( \sigma_1, \dots , \sigma_d) \in \R^d$, and $r$ is the risk-free rate.

The value at time $t$ of an American option with payoff $g$ and maturity $T$ is given by:
\begin{flalign*}
u(t,x) = \sup_{\tau \in {\cal T}_{t,T}} \E[ e^{-r \tau} g(X_{\tau})],
\end{flalign*}
where ${\cal T}_{t,T}$ is the set of stopping time with values in $[t,T]$, and is solution of the variational inequality
\begin{equation}
\left\{
\begin{array}{rcl}
\min \big[ -\partial_t u- \hat \Lc u ,  u -g \big]  &=& 0, \;\;\; \mbox{ on } [0,T)\times (0,\infty)^d \\
 u(T,.) &=& g, \hspace{1cm} \text{ on  } (0,\infty)^d,
\end{array}
\right.
\end{equation}
with 
\begin{align*}
\hat \Lc u(t,x)  & = \;  \frac{1}{2} \sum_{i=1}^d \sigma_i^2 x_i^2 D^2_{x_i} u(t, x)+ r \sum_{i=1}^d x_i D_{x_i} u(t, x) - r u(t,x),    
\end{align*}
as proved e.g. in \cite{jai1990variational}.

Let us define the change of function $v$ by: 
$u(t,x)$ $=$ $e^{ r t} v(t, \log(x))$, (where $\log$ is applied component-wise),  which is solution of the following variational inequality
\begin{equation}
\left\{ 
\begin{array}{rcl}
 \min \left( -\partial_t v-\Lc v ,  v - \hat g \right) &=&0, 
 \;\;\; \mbox{ on } [0,T)\times\R^d \\
\quad v(T,.) &=&\hat g, \hspace{1cm} \text{ on } \R^d,
\end{array}
\right.
 \label{eq:amSimp}
\end{equation}
where $$\hat g(t,x) = e^{-r t} g(e^{x}) ,$$
 $$ \Lc v = \frac{1}{2} \sum_{i=1}^d \sigma_i^2  D^2_{x_i} v_{ii} +  \sum_{i=1}^d  (r- \frac{1}{2}\sigma_i^2)  D_{x_i} v_i.$$
 
In this section, we test the scheme described in section \ref{sec:varIneq} on \eqref{eq:amSimp} in the special case of a geometrical put with strike $K=1$ ,  $T=1$,  $r=0.05$, $X_0^i=1$, $\sigma_i =0.2$ for $i=1$ to $d$, and payoff $(K- \prod_{i=1}^d X^i_t)_{+}$, as considered previously in \cite{bouchard2012monte}. In dimension $d$, the case boils down to the resolution of an American option in dimension $d$ $=$ $1$: indeed, the option payoff involving only the product of the asset values, it can be written as the payoff of a single asset with a trend equal to $d r$ and a volatility $\sigma_1 \sqrt{d}$, so that it can be very accurately estimated e.g. with a tree-based method.
 Results given in Table \ref{tab:American} show that scheme \eqref{eq:schemeVI} is very accurate for the pricing of  American options.
 
 \begin{table}[H]
 \centering
 \begin{tabular}{|c|c|c|c|c|}  \hline
  Dimension &  nb step &  value & std & reference \\ \hline
   1   & 10  &  0.06047  &  0.00023 &   0.060903  \\ \hline
   1   & 20  &  0.060789  & 0.00021   &   0.060903  \\ \hline
   1   & 40  &   0.061122 &  0.00015 &   0.060903  \\ \hline
   1   & 80  &   0.0613818 &  0.00019 &   0.060903  \\ \hline
   5   & 10  &  0.10537 &   0.00014 & 0.10738 \\ \hline
    5   & 20  &  0.10657 &   0.00011& 0.10738 \\ \hline
    5   & 40  &   0.10725&   0.00012 & 0.10738  \\ \hline 
    5   & 80  &   0.107650 &  0.00016 &  0.10738  \\ \hline 
    10  &  10  &  0.12637  &  0.00014   &   0.12996     \\ \hline  
    10  &  20  &  0.128292 &   0.00011  &   0.12996     \\ \hline  
    10  &  40  &  0.12937   &    0.00014       &    0.12996    \\ \hline  
    10  &  80  &  0.129923  &    0.00016     &    0.12996    \\ \hline  
    20 &  10  &   0.1443    &   0.00014    &    0.1510 \\ \hline
    20 &  20  &    0.147781  &   0.00012     &    0.1510 \\ \hline
    20 &  40  &   0.149560 &   0.00012    &    0.1510 \\ \hline
    20 &  80  &   0.15050    &  0.00010    &    0.1510 \\ \hline
    40 & 10 &    0.15512  &  0.00018     &   0.1680  \\ \hline
    40 & 20 &    0.16167  &   0.00015   &   0.1680  \\ \hline
    40 & 40 &    0.16487  &   0.00011   &   0.1680  \\ \hline
    40 & 80 &    0.16665   &   0.00013   &   0.1680  \\ \hline
    40 & 160 &    0.16758  &   0.00016   &   0.1680  \\ \hline

 \end{tabular}
\caption{\label{tab:American} Estimates of the American option using RDBDP. Average and standard deviation over $40$ independent runs for different numbers of time steps are reported.}
\end{table}

%\newpage

\vspace{5mm}

\small

%\bibliography{bsde_references}
%\bibliographystyle{plain}

\printbibliography

\end{document}